%% file: fastSSMpaper.tex
\newif\ifpreprint
\preprinttrue

\ifpreprint
\documentclass[%
 aip,
 cha,
 amsmath,amssymb,
 preprint,%
]{revtex4-1}
\else
\documentclass[%
 aip,
 cha,
 amsmath,amssymb,
 reprint,%
]{revtex4-1}
\fi

\input{packages}
\input{commands}
\usepackage{dcolumn}% Align table columns on decimal point

\usepackage[utf8]{inputenc}
\usepackage[T1]{fontenc}
\usepackage{mathptmx}
\usepackage{etoolbox}

\usepackage{float}

\makeatletter
\def\@email#1#2{%
 \endgroup
 \patchcmd{\titleblock@produce}
  {\frontmatter@RRAPformat}
  {\frontmatter@RRAPformat{\produce@RRAP{*#1\href{mailto:#2}{#2}}}\frontmatter@RRAPformat}
  {}{}
}%
\makeatother
\begin{document}

\preprint{AIP/123-QED}

\title[Fast data-driven model reduction for nonlinear dynamical systems]{Fast data-driven model reduction for nonlinear dynamical systems}
% Force line breaks with \\
\author{Joar Axås}
\author{Mattia Cenedese}%
\author{George Haller}
 \email[Corresponding author: ]{georgehaller@ethz.ch}
\affiliation{ 
Institute for Mechanical Systems, ETH Zürich,
Leonhardstrasse 21, 8092 Zürich, Switzerland
}%

\date{\today}% It is always \today, today,

\begin{abstract}
We present a fast method for nonlinear data-driven model reduction of dynamical systems onto their slowest nonresonant spectral submanifolds (SSMs). 
We use observed data to locate a low-dimensional, attracting slow SSM and compute a maximally sparse approximation to the reduced dynamics on it.
The recently released \SSML{} algorithm uses implicit optimization to fit a spectral submanifold to data and reduce the dynamics to the normal form.
Here, we present two simplified algorithms, which reformulate manifold fitting and normal form computation as explicit problems under certain assumptions.
We show on both numerical and experimental datasets that these algorithms yield accurate and sparse rigorous models for essentially nonlinear (or \textit{non-linearizable}) dynamics.
The new algorithms are significantly simplified and provide a speedup of several orders of magnitude.
\end{abstract}

\maketitle

\ifpreprint
\else
\begin{quotation}
An outstanding challenge in nonlinear dynamics is the development of reduced-order models from data sets representing complex physical systems. 
Here we present a model reduction method that addresses this challenge by constructing a unique, attracting slow manifold (spectral submanifold) for nonlinear systems purely from data. 
The reduced and appropriately sparsified dynamics on this manifold then serves as a mathematically exact reduced-order model for a large set of initial conditions. 
In comparison to prior algorithms targeting spectral submanifolds, the advance here is to use singular value decomposition to compute the tangent space of the manifold and the recursive construction of the sparse reduced dynamics on the manifold via classic normal form theory. 
These simplifications lead to a major speedup which allows for fast and accurate reduced-order modeling from large data sets. 
We illustrate these advantages on experimental data from fluid sloshing and resonant structural vibrations, as well as numerical data obtained from a finite-element beam model.
\end{quotation}
\fi

\section{Introduction}
Nonlinear dynamical systems are omnipresent in nature and engineering. 
Examples include beam and plate buckling, \cite{abramian20} turbulent fluid flows, \cite{holmes12} vibrations in jointed structures, \cite{lacayo19} sloshing in fluid tanks, \cite{taylor53} and even traffic jams. \cite{orosz06}
As computational resources have grown, so has the interest in data-driven methods, which take input data from experiments or simulations and return a reduced model of the underlying system dynamics. 
To date, however, no rigorous method has been accepted as a standard for nonlinear system identification and reduced modeling.

Model simplicity (or parsimony) is vital for interpretability, control, and response prediction for mechanical devices. \cite{kutz22}
This has motivated reduction methods based partially or fully on linearization of the underlying dynamics, such as the proper orthogonal decomposition \cite{lumley67, awrejcewicz04} and the dynamic mode decomposition (DMD). \cite{schmid10, schmid22}
Specifically, DMD obtains the best fit of a linear dynamical system to the data in an equation-free manner, \cite{kutz16} often utilizing delay embedding to secure sufficiently many observables. \cite{dylewsky22}

Many nonlinear mechanical systems, however, exhibit phenomena that cannot occur in linear systems, such as coexisting isolated steady states, and therefore cannot be captured by linear models. \cite{page19}
We refer to such phenomena as \textit{non-linearizable}.
Machine learning methods can potentially capture such phenomena, \cite{brunton16b, lusch18, hartman17, daniel20} but tend to provide models that lack interpretability and perform poorly outside their training range. \cite{loiseau20}

Here, we propose spectral submanifolds (SSMs) to obtain sparse models of nonlinearizable phenomena. 
An SSM is the smoothest nonlinear continuation of a nonresonant spectral subspace emanating from a steady state, both in autonomous systems and systems with periodic or quasiperiodic forcing.
SSMs are unique, attracting, and persistent attracting invariant manifolds in the phase space, that lend themselves well to model reduction. \cite{haller16}
Closely linked to SSMs, invariant spectral foliations are the basis of another rigorous approach for extending linear modal analysis to nonlinear systems. \cite{szalai20}

An available algorithm \cite{SSMTool} computes spectral submanifolds for mechanical systems defined by differential equations. \cite{ponsioen18}
This algorithm takes a Taylor expansion of a set of autonomous ODEs and computes the coefficients of an SSM up to any order and dimension. 
The methodology has been applied to model reduction of systems with hundreds of thousands of degrees of freedom, accurately predicting responses to small harmonic forcing, \cite{ponsioen20, jain18, jain21} as well as the bifurcations of those responses. \cite{li21a, li21b}

These developments have motivated applications of SSM theory to locate reduced-order models from data.
As a first step, Ref~\onlinecite{szalai17} fitted a multivariate polynomial to the sampling map and computed SSMs on the resulting Taylor expanded dynamics.
This yielded good reconstructions of the backbone curves in a clamped-clamped beam experiment.
Typically, however, due to the rapid growth in the number of terms, a polynomial expansion of the flow map in the full observable space is sensitive to overfitting and quickly becomes intractable in higher-dimensional systems.

Recently, the numerical method \SSML{} was developed for data-driven SSM identification. \cite{cenedese21}
This method alleviates the overfitting of the sampling map by separating the model reduction into two steps: manifold identification and reduced dynamics fitting. 
Thus, the polynomial fit of the dynamics takes place only on the SSM, and hence the number of terms no longer grows with the data dimensionality.
First, the data is embedded in a user-defined observable space. 
Then a polynomial representation of the SSM is fitted to the data and the data dimensionality is reduced by projection onto the tangent space of the SSM. 
Using nonlinear optimization techniques, finally a transformation from the reduced coordinates to a normal form is computed, maximizing sparsity while retaining essential nonlinearities. \cite{guckenheimer83}

While data-driven SSM-based model reduction has been successfully applied to both numerical and experimental data for fluid problems, structural dynamics, and fluid-structure interaction, \cite{cenedese21, cenedese21b, kaszas22} the required implicit optimization can be computationally demanding for high-dimensional systems. 
This paper introduces two alternative, simplified algorithms that explicitly compute an SSM followed by a normal form transformation on the reduced coordinates. 
Our key assumption is that the tangent space of the SSM at a fixed point can be found by singular value decomposition (SVD). 
This enables a major simplification and speedup at the expense of giving up the more general applicability of the previous method.

Our new algorithm, \mSSM{}, contains an explicit implementation of the cubic normal form for a complex conjugate eigenvalue pair, suited for monomodal oscillatory dynamics on two-dimensional SSMs. 
Its extension, \mSSMp{}, computes the normal form up to any order and dimensionality.
Both algorithms are fully equivalent in the case of cubic normal forms on 2D SSMs.

Our algorithms compute the normal form on SSMs analytically from a numerical fit of the reduced dynamics. 
This is in contrast to previous methods, which fit the normal form directly to data simultaneously at all required orders.
The aims of our new algorithms are accessibility to practitioners, major speedup for rapid prototyping, analysis of higher-dimensional observable spaces, and insight into the differences between numerically and analytically computed normal forms. 

The structure of this paper is as follows. 
The theory is laid out in \autoref{sec:theory}, with a brief introduction to SSMs, a summary of the \SSML{} algorithm, and a detailed description of \mSSM{} and \mSSMp{}. 
In \autoref{sec:prepost}, we detail how to select an observable space with delay embedding and show how to use SSM-reduced models for forced response prediction. 
In \autoref{sec:examples}, we apply the two newly derived algorithms to experiments from a sloshing tank, simulations of a \vk{} beam, and experiments on an internally resonant beam.
Finally, in \autoref{sec:conclusions} we draw conclusions from the examples, suggest additional applications, and outline possible further enhancements to our dynamics-based machine learning model reduction algorithms.

\section{Model order reduction on spectral submanifolds}\label{sec:theory}
%As shown in \cite{haller16}, under well-defined conditions, spectral subspaces of a linearized system admit unique nonlinear continuations under addition of the higher-order terms. 
%In particular, a \textit{spectral submanifold} (SSM) is the unique smoothest invariant manifold tangent to a spectral subspace. 
%Here we summarize the results for autonomous Taylor-expanded systems.
%
We consider nonlinear autonomous dynamical systems of class $\mathcal{C}^{l}$, $l\in \{\N^+,\ \infty,\ a\}$, with $a$ denoting analyticity, in the form
\begin{equation}\label{eq:fullsystem}
    \dot{\vct{y}} = \vct{A}\vct{y} + \vct{f}(\vct{y}), \quad \vct{y} \in \R^p, \quad \vct{f} \sim \ordo{|\vct{y}|^2}, \quad \vct{f} : \R^p\to \R^p.
\end{equation}

We assume $\vct{A}\in \R^{p\times p}$ is diagonalizable and the real parts of its eigenvalues are all non-zero. 
We denote by $\specsub{}$ the slowest $d$-dimensional spectral subspace of $\vct{A}$, i.e., the span of the eigenvectors corresponding to the $d$ eigenvalues with real part closest to zero. 
When all eigenvalues have negative real parts, the dynamics of the linearized system will converge to $\specsub{}$ in forward time. 

If the eigenvalues corresponding to $\specsub{}$ are non-resonant with the remaining $p-d$ eigenvalues of $\vct{A}$, $\specsub$ admits a unique smoothest, invariant, nonlinear continuation $\mfd$ under addition of the higher-order terms. \cite{cabre03}
We refer to $\mfd$ as a spectral submanifold (SSM). \cite{haller16, cenedese21b}
In the case of a resonance between $\specsub$ and the rest of the spectrum of $\vct A$, we can include the resonant modal subspace in $\specsub$ and thus obtain a higher-dimensional SSM.
If all eigenvalues of $\vct A$ have negative real parts, $\mfd$ will be an attracting slow manifold for system (\ref{eq:fullsystem}), just as $\specsub$ is for the linear part of (\ref{eq:fullsystem}).

A numerical package, \SSMT{}, has been developed for the computation of SSMs from arbitrary finite-dimensional nonlinear systems. \cite{SSMTool, jain21}
More recently, a data-driven algorithm, \SSML{}, has been developed to compute SSMs purely from observables of the dynamical system. \cite{cenedese21, cenedese21b}
In the next section, we first review the \SSML{} algorithm and then discuss two simplified versions of it that enable faster computations under some further assumptions.

\subsection{Learning spectral submanifolds from data: \SSML{}}\label{subsec:ssmlearn}
We seek a parametrization of $\mfd$ from trajectories in an observable space, which may be the full phase space or a suitable embedding, as described in \autoref{subsec:embed}.
All our proposed data-driven methods consist of two steps: manifold fitting and normal form computation.
While \SSML{} solves these problems with implicit optimization, the other two algorithms we develop here perform all these computations explicitly. 

With \SSML{}, \cite{SSMLearn, cenedese21} we compute a parametrization of the SSM as a multivariate polynomial \begin{equation}
\begin{aligned}
    \vct y(\vct \xi) &= \vct M\vct \xi^{1:m} = \vct V\vct \xi + \vct M_{2:m}\vct \xi^{2:m}, \\ 
    \vct{M} &= [\vct{V}, \vct{M}_2, \dots, \vct{M}_m], \quad \vct{M}_i \in \R^{p\times d_i},
\end{aligned}
\end{equation}
with the reduced coordinates $\vct \xi = \vct V^\top \vct y$ and the tangent space $\vct V \in \R^{p\times d}$. 
Here, $d_{i}$ denotes the number of $d$-variate monomials at order $i$.
Throughout this paper, the superscript $(\cdot)^{l:r}$ denotes a vector of all monomials at orders $l$ through $r$. For example, if $\vct{x} = [x_1, x_2]^\top$, then 
$$\vct{x}^{1:3} = [x_1,x_2,x_1^2,x_1x_2,x_2^2,x_1^3,x_1^2x_2,x_1x_2^2,x_2^3]^\top.$$

Minimization of the distance of the parametrized SSM from the training data points $\vct{y}_j$ yields the optimal coefficient matrix $\vct M^\star \in \R^{p\times d_{1:m}}$, with $d_{1:m}$ denoting the number of $d$-variate monomials from orders 1 up to $m$, as
\begin{equation}\label{eq:geometricerror}
\begin{aligned}
    \vct M^\star = [\vct V^\star,\vct M^\star_{2:m}] &= 
    \\
    \mathop{\mathrm{argmin}}_{\vct V,\vct M_{2:m}} \sum_j&||\vct y_j - \vct V\vct V^\top \vct y_j - \vct M_{2:m}(\vct V^\top \vct y_j)^{2:m}||,
\end{aligned}
\end{equation}
subject to the constraints
\begin{equation}
    \vct V^\top \vct V = \vct I, \quad \vct V^\top \vct M_{2:m} = \vct 0.
\end{equation}

Next, we compute the reduced dynamics $\dot{\vct{\xi}}=\vct{R}\vct{\xi}^{1:r}$ on the SSM in polynomial form, with coefficients in the matrix $\vct R\in \R^{d\times d_{1:r}}$, by minimizing
\begin{equation}
    \mathop{\mathrm{argmin}}_{\vct R} \sum_j ||\dot{\vct \xi}_j - \vct R\vct \xi_j^{1:r}||.
\end{equation}

We now compute the normal form of the reduced dynamics of $\vct{\xi}$ on the SSM. \cite{guckenheimer83}
We seek a nonlinear transformation $\vct{t}^{-1}: \vct\xi \mapsto \vct z$ to new coordinates $\vct z \in \Ce^d$, that reduces the number of coefficients in $\vct R$ to a minimum set $\vct N \in \Ce^{d\times d_{1:n}}$.
The transformation and normalized dynamics are given by
\begin{equation}
\begin{aligned}
    \vct z &= \vct t^{-1}(\vct \xi) = \vct H \vct \xi^{1:n} = \vct{W}^{-1} \vct \xi + \vct H_{2:n} \vct \xi^{2:n}, \\
    \dot{\vct z} &= \vct n(\vct z) = \vct N\vct z^{1:n} = \vct \Lambda \vct z + \vct N_{2:n}\vct z^{2:n},
\end{aligned}
\end{equation}
where $\vct{H}\in \Ce^{d\times d_{1:n}}$, and $\vct W \in \Ce^{d\times d}$ is the matrix of eigenvectors of the linear part $\vct R_1 = \vct W \vct \Lambda \vct{W}^{-1}$ of the reduced dynamics.
$\vct N$ is potentially a very sparse matrix, simplifying the dynamics $\vct n(\vct z)$.
The locations of nonzero elements in $\vct N$ are determined by any approximate inner resonances between the eigenvalues in $\vct \Lambda$. \cite{cenedese21}
An example is given in (\ref{eq:Nforcubic}).

In \SSML{}, the transformation to the normal form is computed by minimization of the conjugacy error
\begin{equation}\label{eq:conjugacyerror}
    \mathop{\mathrm{argmin}}_{\vct N,\vct H} \sum_j||\vct \nabla_{\vct \xi} \vct t^{-1}(\vct \xi_j)\dot{\vct \xi}_j - \vct n(\vct t^{-1}(\vct \xi_j))||.
\end{equation}

\subsection{Fast SSM identification -- \mSSM{} and \mSSMp{}}\label{subsec:minissm}
\begin{figure*}
    \centering
    \includegraphics[width=\textwidth]{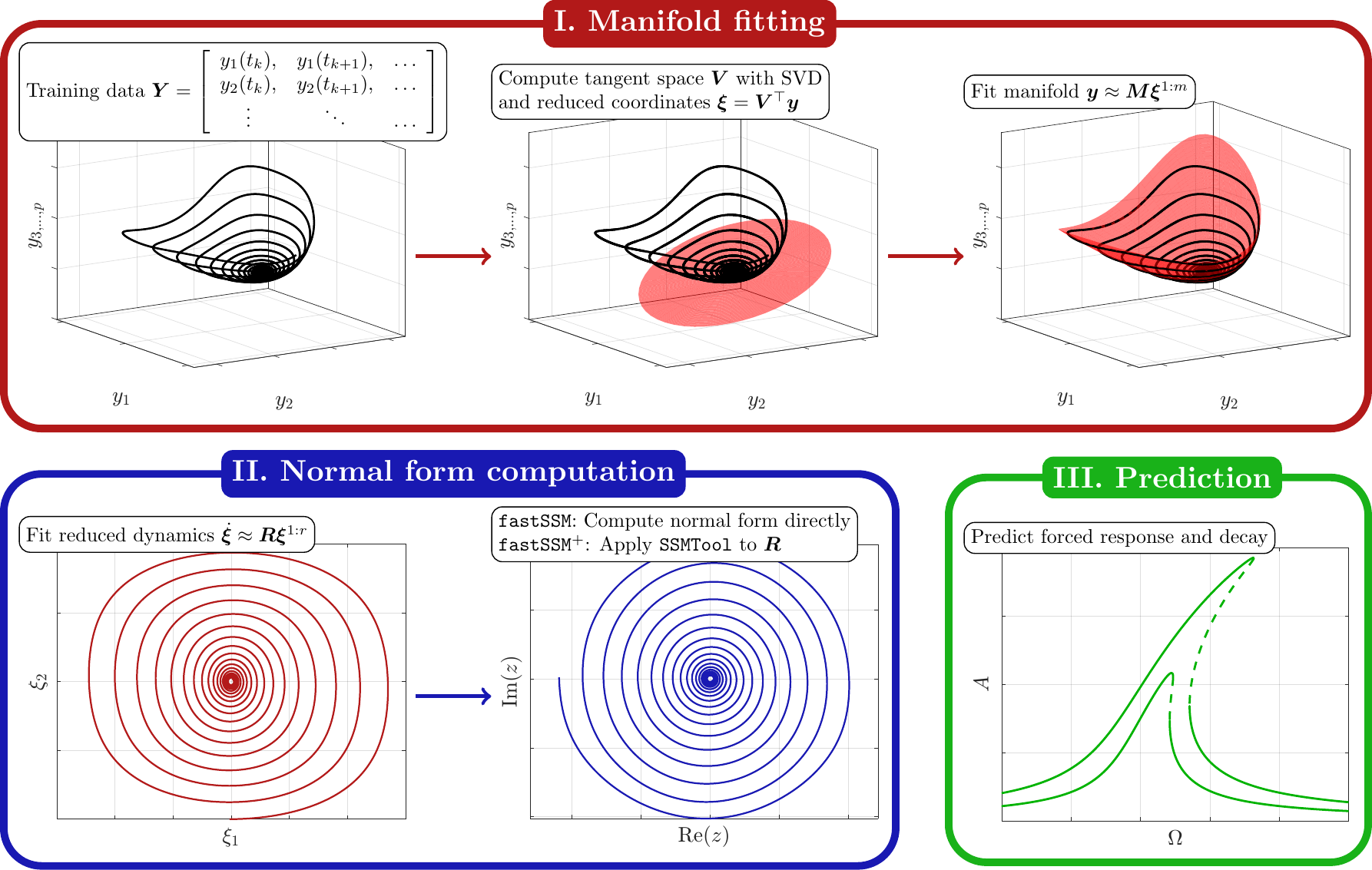}
    \caption{The \mSSM{} algorithm consists of two main parts, (I) fitting an SSM with SVD and polynomial regression, and (II) computing a normal form on the manifold using \SSMT{}. 
    The obtained reduced model can then be used for (III) prediction of the autonomous system evolution or the response of the forced system.}
    \label{fig:mssm}
\end{figure*}

To turn the manifold fitting into an explicit problem, we assume that the tangent space of the SSM at the origin can be approximately obtained by singular value decomposition (SVD) on the data. 
This is typically satisfied sufficiently close to an equilibrium point.
Moreover, since data variance prevalently tends to occur along the tangent space of the SSM, SVD is a good candidate for this subspace.
A third motivation is that the image of $\mfd$ in a delay-embedded space tends to be flat, even if $\mfd$ is strongly curved in the full phase space. \cite{cenedese21}
Although the computed manifold may curve significantly, it is not allowed to fold over its tangent space. 
An overview of the algorithm is shown in \autoref{fig:mssm}. 

The tangent space $\specsub{}$ is approximated by $d$-dimensional truncated SVD on the snapshot matrix $\vct Y \in \R^{p\times N}$ as 
\begin{equation}
    \vct Y \approx \vct U_d\vct S_d\vct{\hat V}_d^\top,
\end{equation}
where $\vct U_d \in \R^{p\times d}$, $\vct S_d \in \R^{d\times d}$, $\vct{\hat V}_d \in \R^{N\times d}$, and by defining the tangent space approximation as
\begin{equation}
    \vct V = \vct U_d\vct S_d^{-1},
\end{equation}
we can project onto it to obtain $d$ reduced coordinates $\vct \xi$ as
\begin{equation}
    \vct \xi = \vct V^\top \vct y.
\end{equation}
We denote by $\vct \Xi \in \R^{d\times N}$ the projection of $\vct Y$ onto $\vct V$. 
Numerically, it is beneficial to first normalize each column of $\vct V$ with the maximum absolute value of the corresponding row in $\vct{\hat V}_d$

The manifold parametrization coefficients $\vct M\in \R^{p\times d_{1:m}}$ are obtained by polynomial regression: 
\begin{equation}\label{eq:geometric}
    \vct y \approx \vct M\vct \xi^{1:m}, \quad \vct M = \vct Y(\vct \Xi^{1:m})^\dagger,
\end{equation}
with $(\cdot)^\dagger$ denoting the Moore-Penrose pseudo-inverse.

The reduced dynamics are approximated by an $\ordo{r}$ polynomial regression, with the coefficient matrix $\vct R \in \R^{d\times d_{1:r}}$, as
\begin{equation}\label{eq:reddyn}
    \dot{\vct \xi} \approx \vct R\vct \xi^{1:r}, \quad \vct R = \dot{\vct \Xi}(\vct \Xi^{1:r})^\dagger.
\end{equation}
The time derivative $\dot{\vct \xi}$ may be computed numerically. 
The quality of this approximation is important for model accuracy. 
Therefore, we use a central finite difference method including 4 adjacent points in each time direction. \cite{fornberg88}

Next, diagonalizing the linear part $\vct R_1 = \vct W\vct \Lambda \vct W^{-1}$, we again apply regression to rewrite $\vct R$ in the form $\vct G\in \Ce^{d\times d_{1:r}}$ with the eigenvalues in the linear part $\vct G_1=\vct \Lambda$ and modal coordinates $\vct \zeta = \vct{W}^{-1}\vct z$, such that 
\begin{equation}\label{eq:diagonalize}
    \dot{\vct \zeta} \approx \vct G\vct{\zeta}^{1:r}.
\end{equation}

To compute the normal form of (\ref{eq:diagonalize}), we seek a near-identity polynomial transformation $\vct t: \vct z \mapsto \vct \zeta$ with coefficients $\vct T \in \Ce^{d\times d_{1:n}}$ from new coordinates $\vct z \in \Ce^d$ such that
\begin{equation}\label{eq:generalnf}
\begin{aligned}
     \vct \zeta &= \vct t(\vct z) = \vct T\vct z^{1:n} = \vct z + \vct T_{2:n}\vct z^{2:n} \\
    \dot{\vct z} &= \vct n(\vct z) = \vct N\vct z^{1:n} = \vct \Lambda \vct z + \vct N_{2:n}\vct z^{2:n}.
\end{aligned}
\end{equation}
We enforce conjugacy between the normal form and reduced dynamics by plugging (\ref{eq:generalnf}) into (\ref{eq:diagonalize}):
\begin{equation}\label{eq:conjugacy}
    \vct \nabla_{\vct z}(\vct T\vct z^{1:n})\vct N\vct z^{1:n} = \vct G(\vct T\vct z^{1:n})^{1:r}.
\end{equation}
We can compute $\vct T$ and $\vct N$ by solving (\ref{eq:conjugacy}) recursively at increasing orders. 
See Ref.~\onlinecite{jain21} for details on the computations.

The simplest non-trivial normal form arises on a 2D manifold emanating from a spectral subspace corresponding to two complex conjugate eigenvalues $(\lambda,\lbar)$ with small real part. 
We denote the coordinates in the normal form $(z_1,z_2)=(z,\zbar)$. 
The cubic normal form, with $\gamma\in \Ce$, is
\begin{equation}\label{eq:cubicnf}
\begin{aligned}
    \dot{z} &= \lambda z + \gamma z^2\zbar + \ordo{5} \\
    \dot{\zbar} &= \lbar \zbar + \bar{\gamma} z\zbar^2 + \ordo{5},
\end{aligned}
\end{equation}
for which we obtain
\begin{equation}\label{eq:Nforcubic}
    \vct N = \left[\begin{array}{ccccccccc}\lambda & 0 & 0 & 0 & 0 & 0 & \gamma & 0 & 0\\
     0 & \lbar & 0 & 0 & 0 & 0 & 0 & \bar{\gamma} & 0\end{array}\right].
\end{equation}
Solving (\ref{eq:conjugacy}) yields
\begin{equation}
    \vct T = \left[\begin{array}{ccccccccc}1 & 0 & T_{3} & T_{4} & T_{5} & T_{6} & 0 & T_{8} & T_{9}\\
     0 & 1 & \bar{T}_{5} & \bar{T}_{4} & \bar{T}_{3} & \bar{T}_{9} & \bar{T}_{8} & 0 & \bar{T}_{6}\end{array}\right],
\end{equation}
where
\begin{equation}
\begin{aligned}
    T_3 &= \frac{G_3}{\lambda}, \quad T_4 = \frac{G_4}{\lbar}, \quad T_5 = \frac{G_5}{2\lbar-\lambda}, \\
    T_6 &= \frac{2G_3 T_3 + G_4 \bar{T}_5 + G_6}{2\lambda}, \quad
    T_9 = \frac{G_4 T_5 + 2G_5 \bar{T}_3 + G_9}{3\lbar-\lambda}, \\
    T_8 &= \frac{2G_3 T_5 + G_4\bar{T}_3 + G_4 T_4 + 2G_5 \bar{T}_4 + G_8}{2\lbar}, \\
    \gamma &= 2G_3 T_4 + G_4 \bar{T}_4 + G_4 T_3 + 2G_5 \bar{T}_5 + G_7,
\end{aligned}
\end{equation}
and $G_i$ refers to element $i$ of the top row in $\vct G$. \cite{szalai17}

This solution is implemented in \mSSM{} for fast, cubic-order, 2D normal form computations. 
\ifpreprint
\else
The full script is written out in Appendix \ref{appendix:minissm}.
\fi
We have the reduced dynamics and normal form orders $r=3$ and $n=3$, and we are free to set the order $m$ of the manifold. 

While a cubic order normal form will suffice to model a number of datasets accurately, \cite{cenedese21, cenedese21b} for stronger nonlinearities we must include higher orders.
\SSMT{} can solve (\ref{eq:conjugacy}) for any dimension and order of expansion.
The algorithm \mSSMp{} extends the manifold fitting and normal form computation in \mSSM{} to any dimension and order by automatically applying \SSMT{} to the coefficients of $\vct R$.
\ifpreprint
\else
The full code is written out in Appendix \ref{appendix:minissmplus}.
\fi
In principle, we are free to choose any orders $m$, $r$, and $n$, but to avoid overfitting we must limit the manifold order $m$ and reduced dynamics order $r$. 
In addition, we need to pick a large enough $n$ to make the tail of the SSM-reduced dynamics small enough.

Unlike (\ref{eq:geometricerror}) and (\ref{eq:conjugacyerror}) in \SSML{}, all computations in \mSSM{} and \mSSMp{} are explicit. 
Therefore, we expect these algorithms to be faster than \SSML{}, at the cost of some model accuracy.
Further, we note that in \mSSM{} and \mSSMp{}, the normal form is computed analytically from a numerical polynomial fit, whereas \SSML{} directly fits the normal form coefficients numerically from data. 

In practice we must also approximate the inverse of $\vct t$ to transform initial conditions to the normal form before the integration. 
The most accurate and simplest way is a numerical solution for each initial condition, but a general inverse map can also be computed.
Accordingly, in \mSSM{}, we implement explicitly the 3rd-order polynomial expansion of the inverse of $\vct t$, whereas \mSSMp{} fits a polynomial map by regression on the training data.

Finally, we note that all our methods require that the data lies sufficiently close to an invariant manifold.
We achieve this by removing any initial transients from the training dataset, as identified by a spectral analysis on the input signal. \cite{cenedese21b}
Since the SSM is unique and attracting, this procedure ensures that we train on relevant data.

\section{Pre- and post-processing of the data}\label{sec:prepost}
We present here considerations for the delay embedding of observable functions to facilitate manifold identification. 
We also show how to use the obtained normal form dynamics to predict the forced response of the full system.

\subsection{Embedding spectral submanifolds in the observable space}\label{subsec:embed}
In practice, observing trajectories $\vct x(t)$ in the full phase space is often intractable.
For instance, an experimentalist might only observe a scalar quantity $s(t) = \zeta(\vct x(t))$. 
Nevertheless, if the observable function $\zeta: \R^{n_\mathrm{full}} \to \R$ is generic, it can be used to reconstruct an invariant manifold of the system in an observable space via delay embedding. 

To this end, we collect $p$ measurements separated by a time lag $\tau>0$ to form a vector $\vct y$ in an observable space $\R^p$ as 
$$\vct y(t) = [s(t), s(t+\tau), s(t+2\tau), \dots, s(t+(p-1)\tau)]^\top.$$ 
Takens' embedding theorem implies that if $\vct x(t)$ lies on an invariant $d$-dimensional manifold $\mfd$ in the original phase space, then $\vct y(t)$ lies on a diffeomorphic copy of $\mfd$ in $\R^p$ with probability one, if $p\ge 2d+1$ and $\zeta(\vct{x})$ is a generic observable. \cite{takens81}
This result also extends to generic multivariate functions $\vct \zeta: \R^{n_\mathrm{full}} \to \R^{n_\mathrm{obs}}$ as long as the total observable space dimension exceeds $2d$. \cite{sauer91, deyle11}
Even with a high-dimensional observable, however, we apply delay embedding to diversify the data and unveil information about the time derivative of the signal in our later sloshing example. \cite{cenedese21b}

The timelag $\tau$ is typically a multiple of the sampling time or timestep $\Delta t$ of the available dataset.
The default choice is $\tau=\Delta t$, but one can also increase the time lag to $\tau=k\Delta t$, $k\in \N_+$. 
Delay-embedding an observed time series $[s(t_1), \dots, s(t_N)]$ thus yields the snapshot matrix
\begin{equation}
    \vct Y = \left[\begin{array}{cccc}
    s(t_1) & s(t_2) & \dots & s(t_{N+(1-p)k}) \\
    s(t_{k+1}) & s(t_{k+2}) & \dots & s(t_{N+(2-p)k}) \\
    s(t_{2k+1}) & s(t_{2k+2}) & \dots & s(t_{N+(3-p)k}) \\
    \vdots & \vdots & \ddots & \vdots \\
    s(t_{(p-1)k+1}) & s(t_{(p-1)k+2}) & \dots & s(t_{N})
    \end{array}\right].
\end{equation}
Takens' embedding theorem requires at least $k=1$, $p=2d+1$, which typically suffices for SSM identification when $d=2$.

If we observe a smooth scalar signal $s(t)$, then for small $\tau$, we have
$$s(t+j\tau) = s(t) + \dot{s}(t)j\tau + \ordo{(j\tau)^2}.$$
Delay-embedding $s(t)$ then yields
\begin{equation}\label{eq:delayplane}
\begin{aligned}
    \vct y(t) = \left[\begin{array}{c}
	s(t)\\
	s(t+\tau)\\
	s(t+2\tau)\\
	\vdots \\
	s(t+(p-1)\tau)
    \end{array}\right]
	=
    s(t)\left[\begin{array}{c}
	1 \\
	1 \\
	1 \\
	\vdots \\
	1 
    \end{array}\right]
	&+
	\dot s(t)\tau\left[\begin{array}{c}
	0\\
	1\\
	2\\
	\vdots \\
	p-1
    \end{array}\right] \\
	&+ \ordo{(p\tau)^2}.
\end{aligned}
\end{equation}
As noted in Ref.~\onlinecite{cenedese21}, if $\tau$ and $p$ are small in comparison to $\dot{s}$, then the diffeomorphic copy of $\mfd$ in the observable space is approximately a plane spanned by the two constant vectors appearing in (\ref{eq:delayplane}). 

For multimodal signals, however, $p$ and $k$ must be increased to allow distinction between the modes in the observable space. 
%Consider a univariate oscillatory signal $s(t) = \zeta(\vct x(t))$ dominated by two modes of oscillation with instantaneous frequencies $\omega_1(t)$,~$\omega_2(t)$ sampled with a small timestep. 
The modal subspaces of the full phase space have corresponding planes in the observable space, which must be identified to provide reduced coordinates for the SSM. 
A low delay embedding dimension and an overly short timelag can make these planes close to parallel, which complicates their identification. 
In contrast, we want to pick the timelag and embedding dimension such that the images of the modal subspaces in the observable space are close to orthogonal. 

To illustrate this, we consider an observed superposition of two harmonics $s(t) = a\cos(\omega_1 t+\psi_1) + b\cos(\omega_2 t+\psi_2)$, a reasonable model of a typical signal close to the origin.
Using a trigonometric identity, we can rewrite
\begin{equation}
\begin{aligned}
	s(t+j\tau) &= 
	a\cos(\omega_1 t+\psi_1)\cos(\omega_1 j\tau) -\\
	 &a\sin(\omega_1 t+\psi_1)\sin(\omega_1 j\tau) +\\ 
	&b\cos(\omega_2 t+\psi_2)\cos(\omega_2 j\tau) -\\
	 &b\sin(\omega_2 t+\psi_2)\sin(\omega_2 j\tau).
\end{aligned}
\end{equation}
The delay-embedded observable vector is therefore
\begin{equation}
    \vct y(t) = 
    \left[\begin{array}{c}
	s(t)\\
	s(t+\tau)\\
	s(t+2\tau)\\
	\vdots \\
	s(t+(p-1)\tau)
    \end{array}\right]
	=
    \vct V
    \left[\begin{array}{c}
    a\cos(\omega_1 t+\psi_1) \\
    -a\sin(\omega_1 t+\psi_1) \\
    b\cos(\omega_2 t+\psi_2) \\
    -b\sin(\omega_2 t+\psi_2)
    \end{array}\right],
\end{equation}
where
\begin{equation}
    \vct V^\top = 
    \left[\begin{array}{ccccc}
	1 & \cos(\omega_1 \tau) & \cos(2\omega_1 \tau) & \dots & \cos((p-1)\omega_1 \tau)\\
	0 & \sin(\omega_1 \tau) & \sin(2\omega_1 \tau) & \dots & \sin((p-1)\omega_1 \tau) \\
	1 & \cos(\omega_2 \tau) & \cos(2\omega_2 \tau) & \dots & \cos((p-1)\omega_2 \tau)  \\
	0 & \sin(\omega_2 \tau) & \sin(2\omega_2 \tau) & \dots & \sin((p-1)\omega_2 \tau)  \\
    \end{array}\right].
\end{equation}
%\begin{widetext}
%\begin{equation}
 %   \vct V = 
  %  \left[\begin{array}{cccc}
	%1 & 0 & 1 & 0\\
	%\cos(\omega_1 \tau) & \sin(\omega_1 \tau) & \cos(\omega_2 \tau) & \sin(\omega_2 \tau) \\
	%\cos(2\omega_1 \tau) & \sin(2\omega_1 \tau) & \cos(2\omega_2 \tau) & \sin(2\omega_2 \tau)  \\
	%\vdots & \vdots & \vdots & \vdots \\
	%\cos((p-1)\omega_1 \tau) & \sin((p-1)\omega_1 \tau) & \cos((p-1)\omega_2 \tau) & \sin((p-1)\omega_2 \tau)  \\
    %\end{array}\right],
%\end{equation}
%\end{widetext}
The trajectories will reside on a 4-dimensional hyperplane given by the \textit{constant} matrix $\vct V$ of delay-embedded harmonic functions.
Picking $\tau$ and $p$ too small complicates manifold identification, because the columns of $\vct V$ are close to linearly dependent. 
Instead, we choose a timelag such that these planes are as close to orthogonal as possible.
Based on these considerations, we select a timelag $\tau = k\Delta t$ such that $\omega_2 \tau \approx \frac{\pi}{2}$, where $\omega_2$ is the second eigenfrequency of the observed system. 
This argument generalizes to any number of modes, although the choice of $k$ becomes more complicated.
The embedding dimension $p$ is set to at least $2d+1$.
Further increases to $p$ can facilitate identification of a manifold but also complicate its geometry far from the origin.
%Although the trajectory will deviate from $V$ under addition of nonlinearities, this method provides an initial guess for $k$.

%Let us assume that the dynamics is dominated by $d/2$ oscillatory modes with coordinates $(z_\ell, \zbar_\ell)$ in the normal form, $\ell = 1,\dots,d/2$. We are observing a scalar function $s(t) = \zeta(\vct m(\vct t(z_\ell(t))))$. Stacking $p$ measurements of $s(t)$ with timelag $\tau = k\Delta t$ yields the observable vector $\vct y(t)$ where the $(j+1)$th element is
%\begin{equation}\label{eq:defy}
%	y_{j+1}(t) = s(t + j\tau) = s(z_l(t + j\tau) = s(\rho_l(t+j\tau)\e^{i\theta_l(t+j\tau)}), \quad j = 0,\dots, (p-1).
%\end{equation}
%We now assume that the derivatives of $\rho_\ell$ and $\theta_\ell$ are approximately constant over any time window $p\tau$ and write $\rho_\ell(t+j\tau) \approx \rho_\ell(t) + \dot{\rho}_\ell(t)j\tau$, $\theta_\ell(t+j\tau) \approx \theta_\ell(t) + \dot{\theta}_\ell(t)j\tau$. 
%Plugging the approximations into (\ref{eq:defy}) yields
%\begin{equation}
%	y_{j+1}(t) \approx 
%	s(\rho_\ell(t)\e^{\theta_\ell(t)}\left(1+\frac{\dot{\rho}_\ell(t)}{\rho_\ell(t)} j\tau\right)\e^{\dot\theta_\ell(t) j\tau}) = 
%	s(z_\ell(t)(1+c_\ell j\tau)\e^{\omega_\ell j\tau}).
%\end{equation}
%In other words, we can factorize the dependencies of $t$ and $\tau$ and trajectories will be confined to a hyperplane given by $V_{j+1,\ell} = (1+c_\ell j\tau)\e^{\omega_\ell j\tau}$

\subsection{Backbone curves and forced response}\label{subsec:frc}
If all eigenvalues of the reduced dynamics are complex conjugate pairs, we can rewrite the normal form on the SSM in polar coordinates $(\vct \rho, \vct \theta)$ with $z_\ell = \rho_\ell\e^{i\theta_\ell}$, $z_{\ell+1} = \rho_\ell\e^{-i\theta_\ell}$. 
This yields the SSM-reduced dynamics in the form
\begin{equation}
\begin{aligned}
    \dot{\rho}_\ell &= \damp_\ell(\vct \rho, \vct \theta)\rho_\ell, \\
    \dot{\theta}_\ell &= \freq_\ell(\vct \rho,\vct \theta).
\end{aligned}
\end{equation}
If the imaginary parts of all eigenvalue pairs are non-resonant with each other, then the functions $\damp_\ell$ and $\freq_\ell$ depend only on $\vct \rho$. \cite{ponsioen18}
For example, the cubic 2D normal form (\ref{eq:cubicnf}) in polar coordinates reads
\begin{equation}
\begin{aligned}
    \dot{\rho} &= \real(\lambda)\rho + \real(\gamma) \rho^3 + \ordo{\rho^5}, \\
    \dot{\theta} &= \imag(\lambda) + \imag(\gamma) \rho^2 + \ordo{\rho^4}.
\end{aligned}
\end{equation}

For the general 2D case, we define an amplitude as
\begin{equation}
    \mathcal{A}(\rho) = \max_{\theta\in [0,2\pi)}\left|\alpha\left(\vct W\vct t\left(\left[\begin{array}{c}\rho\e^{i\theta}\\ \rho\e^{-i\theta}\end{array}\right]\right)\right)\right|,
\end{equation}
where $\alpha: \R^p \to \R$ is a function mapping from the observable space to some amplitude of a particular degree of freedom or the norm of total displacements. 
The backbone curve of that amplitude is then defined as the parametrized curve
\begin{equation}\label{eq:backbone}
    \mathcal{B}(\rho) = (\freq(\rho), \mathcal{A}(\rho)),
\end{equation}
which is broadly used in the field of nonlinear vibrations to illustrate the overall effect of nonlinearities in the system.

Next, following Refs.~\onlinecite{breunung18, ponsioen19}, we use the data-driven normal form on the SSM to predict the response of the system under additional, time-periodic external forcing. This amounts to adding a forcing term with amplitude $f$ and frequency $\Omega$ to the general 2D normal form to obtain
\begin{equation}\label{eq:forcednf}
\begin{aligned}
    \dot{\rho} &= \damp(\rho)\rho + f\sin\psi, \\
    \dot{\psi} &= \freq(\rho) - \Omega + \frac{f}{\rho}\cos\psi,
\end{aligned}
\end{equation}
where we have introduced the phase difference $\psi = \theta - \Omega t$. 
The forced response curve (FRC) is then defined as the bifurcation curve of the fixed points $(\rho_0,\psi_0)$ of (\ref{eq:forcednf}) under varying $\Omega$. 
Squaring and adding the equations in (\ref{eq:forcednf}) yields a parametrization of the FRC for the forcing frequency $\Omega$ and the phase lag $\psi_0$ in the form
\begin{equation}\label{eq:frc}
\begin{aligned}
    \Omega(\rho_0,f) &= \freq(\rho_0) \pm \sqrt{\frac{f^2}{\rho_0^2}-(\damp(\rho_0))^2}, \\
    \psi_0(\rho_0,f) &= -\arcsin\left(\frac{\damp(\rho_0)\rho_0}{f} \right).
\end{aligned}
\end{equation}

%Following \cite{cenedese21, cenedese21b}, we may use the first equation in (\ref{eq:frc}) to calibrate the normal form forcing amplitude $f$ to some true observation. 
%In practice, the mapping of a load in the full phase space to the normal form is unknown, because the model contains no information on how the full system responds to forcing. 
%Therefore, one calibration point $(\vct y_\cal,\Omega)$ per forcing level is needed. 
%The observed response $\vct y_\cal$ can be projected onto the manifold and transformed to the normal form to obtain the calibration normal form amplitude $\rho_\cal$.
%The forced response curve is parametrized by the normal form amplitude $\rho$ and forcing $f$. 
Since the relation between the experimental forcing and the normal form forcing $f$ is unknown, a calibration of each FRC to at least one observation of a forced response is necessary. 
Following Refs.~\onlinecite{cenedese21, cenedese21b}, we achieve this by prescribing a single intersection point on the FRC in the frequency-amplitude plane. 
We use the calibration frequency $\Omega_\cal$ and amplitude $u_\cal$ of the maximal experimental response at each forcing level. 
The mapping to an equivalent point $\vct y_\cal$ in the observable space can then be approximated by delay-embedding a cosine signal of amplitude $u_\cal$ and frequency $\Omega_\cal$. 
We compute the calibration amplitude $\rho_\cal$ by projecting $\vct y_\cal$ onto the manifold and transforming to the normal form. 
The forcing $f$ is then computed from the relationship
\begin{equation}
	f^2 = (\damp(\rho_\cal))^2\rho_\cal^2 + (\freq(\rho_\cal) - \Omega)^2\rho_\cal^2.
\end{equation}

Finally, under periodic forcing, the SSM parametrization becomes time-dependent with the addition of a small periodic term. 
We ignore this contribution here to simplify the analysis, but note that this small non-autonomous correction can be exactly computed using \SSML{}. \cite{cenedese21, cenedese21b}

\section{Examples}\label{sec:examples}
We now apply \mSSM{} and \mSSMp{} to three datasets: sloshing experiments in a water tank, simulation of a clamped-clamped \vk{} beam, and experiments on an internally resonant beam. 

We quantify the quality of our SSM-reduced models with the \textit{normalized mean trajectory error} (NMTE). \cite{cenedese21} 
Given a test trajectory with $N$ snapshots $\vct y(t_j)$ and the model-based reconstruction $\hat{\vct y}(t_j)$ obtained by integrating the normal form dynamics and mapping back to the observable space, the NMTE is defined as
\begin{equation}
	\mathrm{NMTE} = \frac{1}{N \max_j||\vct y(t_j)||}\sum_{j=1}^N ||\hat{\vct y}(t_j)-\vct y(t_j)||.
\end{equation}

\subsection{Tank sloshing}\label{subsec:sloshing}
\begin{figure}
     \centering
%     \begin{subfigure}{1\linewidth}
         \includegraphics[width=\linewidth]{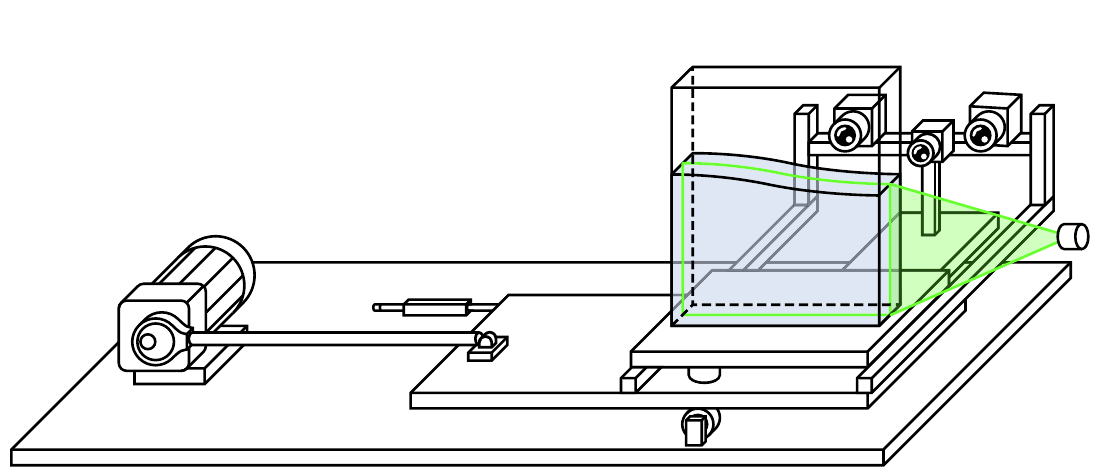}
%         \caption{}
%     \end{subfigure}
%     \begin{subfigure}{0.8\linewidth}
%         \includegraphics[width=\linewidth]{figs/sloshingMiniRec.pdf}
%         \caption{}
%         \label{fig:decay}
%     \end{subfigure}
     \caption{
%     (\subref{fig:experimentsetup}) 
     Sloshing experiment setup with the partially filled tank and the cameras, mounted on a platform excited by a motor. 
%     (\subref{fig:decay}) The decay of the wall surface amplitude agrees with the reconstruction from the reduced dynamics obtained with \mSSM{}. 
     }
%     \label{fig:sloshing1}
         \label{fig:experimentsetup}
\end{figure}

A tank partially filled with a liquid responds nonlinearly to horizontal harmonic excitation. \cite{taylor53}
Stronger fluid oscillation gives rise to more shearing against the tank wall, so that the damping of the system increases nonlinearly with the amplitude. \cite{Faltinsen09}
In addition, the instantaneous frequency decreases at higher amplitudes. 
Both phenomena are crucial for predicting the forced response amplitude of the liquid. 
The industrial applications for sloshing models are numerous, ranging from the transportation of fluids in trucks \cite{cheli13} and container ships \cite{Faltinsen09} to the design of fuel tanks for spacecraft. \cite{dodge00, abramson66}

\begin{figure*}
     \centering
     \parbox{0.32\linewidth}{
     \begin{subfigure}{1\linewidth}
         \includegraphics[width=\linewidth]{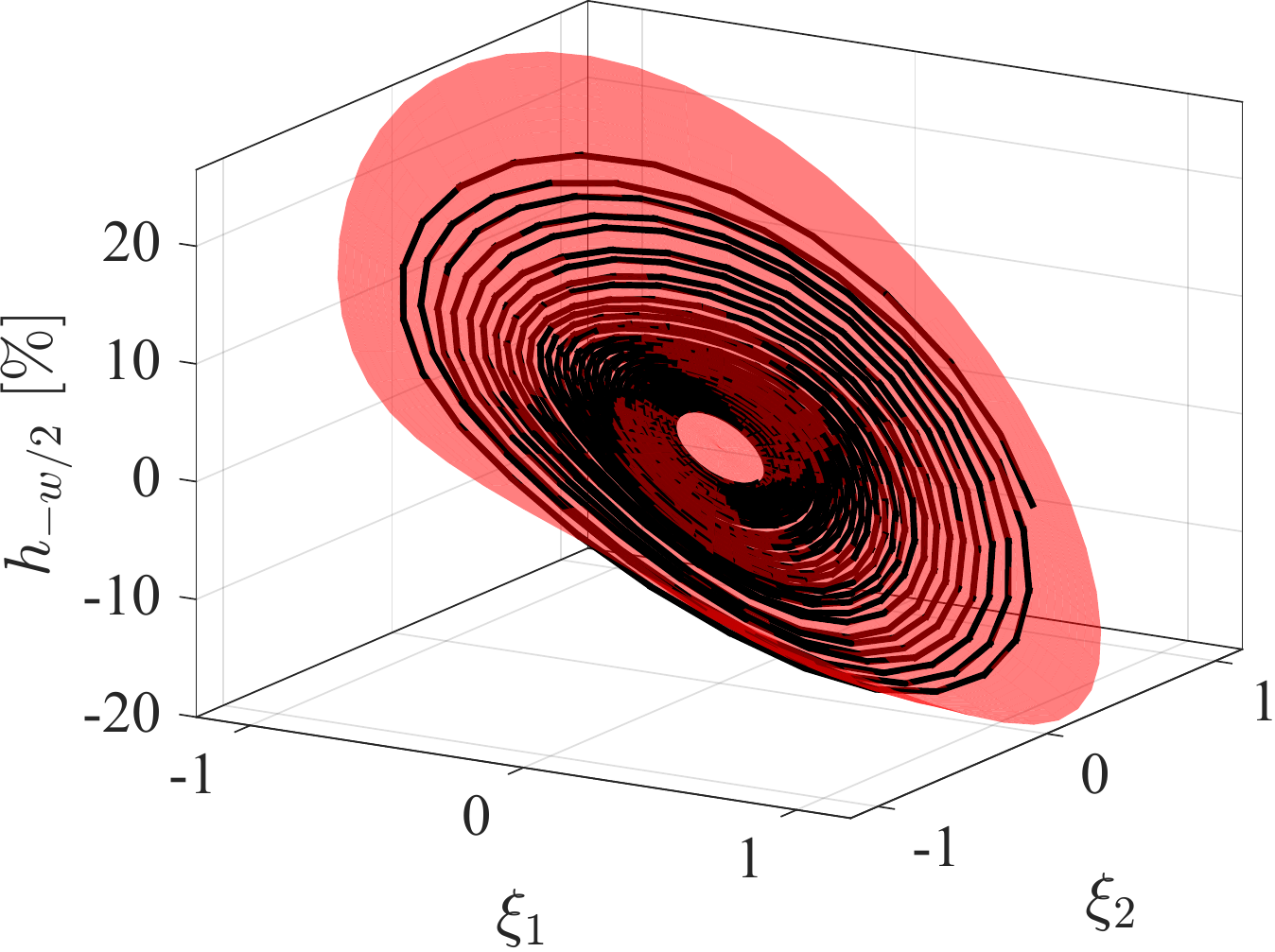}
         \caption{}
         \label{fig:ssm}
     \end{subfigure}
     \begin{subfigure}{1\linewidth}
         \includegraphics[width=\linewidth]{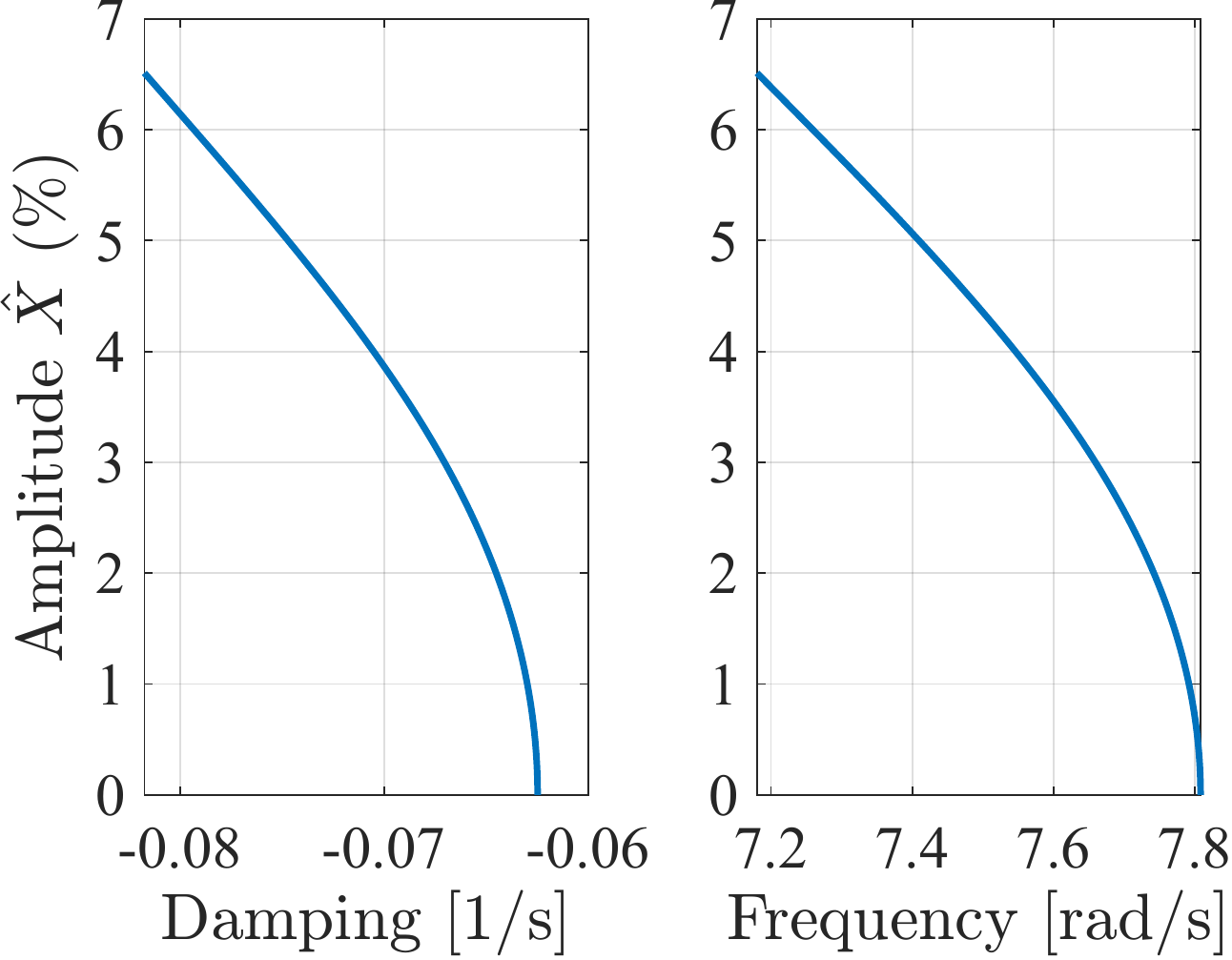}
         \subcaption{}
         \label{fig:backbones}
     \end{subfigure}}
     \parbox{0.32\linewidth}{
     \begin{subfigure}{1\linewidth}
         \includegraphics[width=\linewidth]{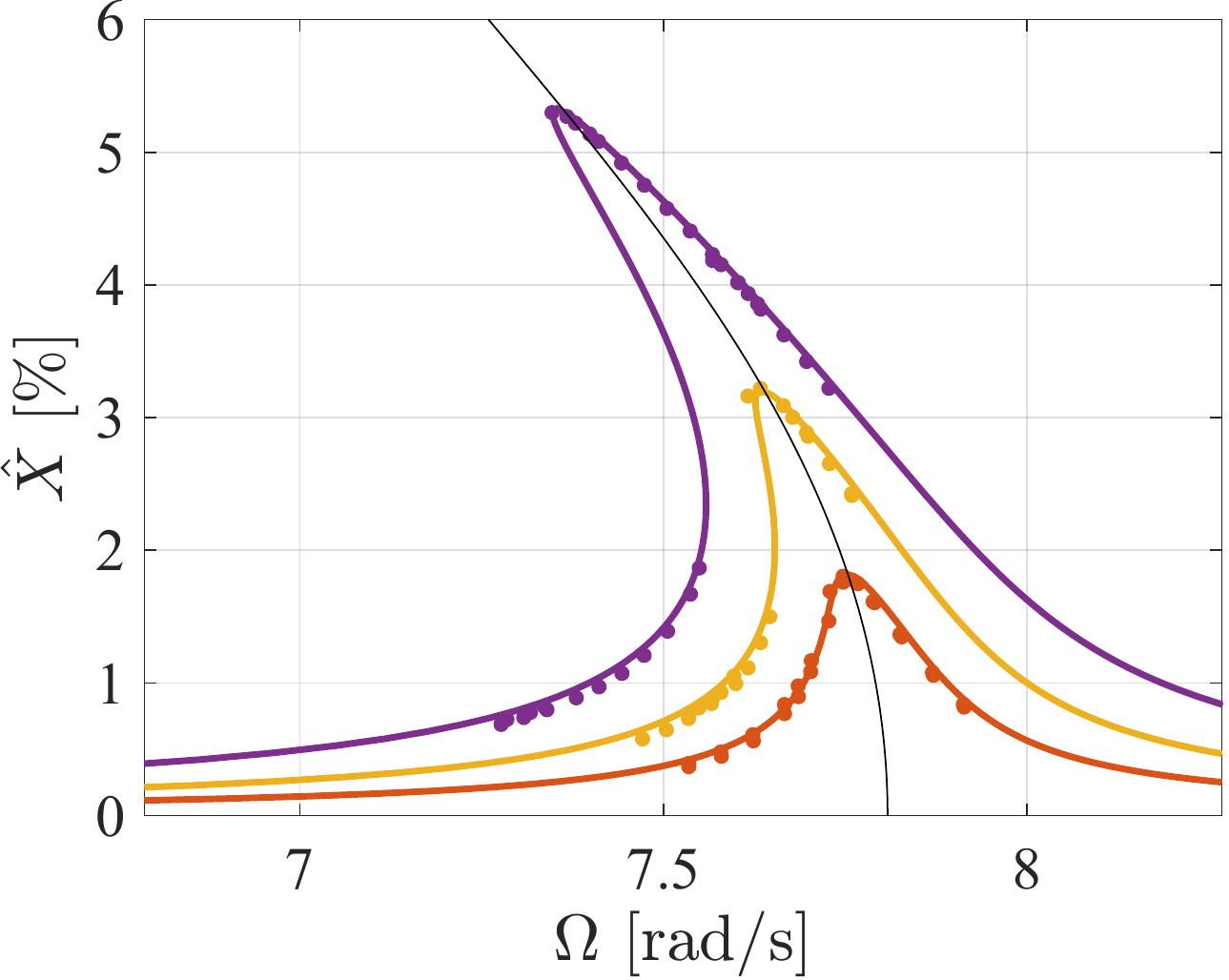}
         \caption{}
         \label{fig:frcamp}
     \end{subfigure}
     \begin{subfigure}{1\linewidth}
         \includegraphics[width=\linewidth]{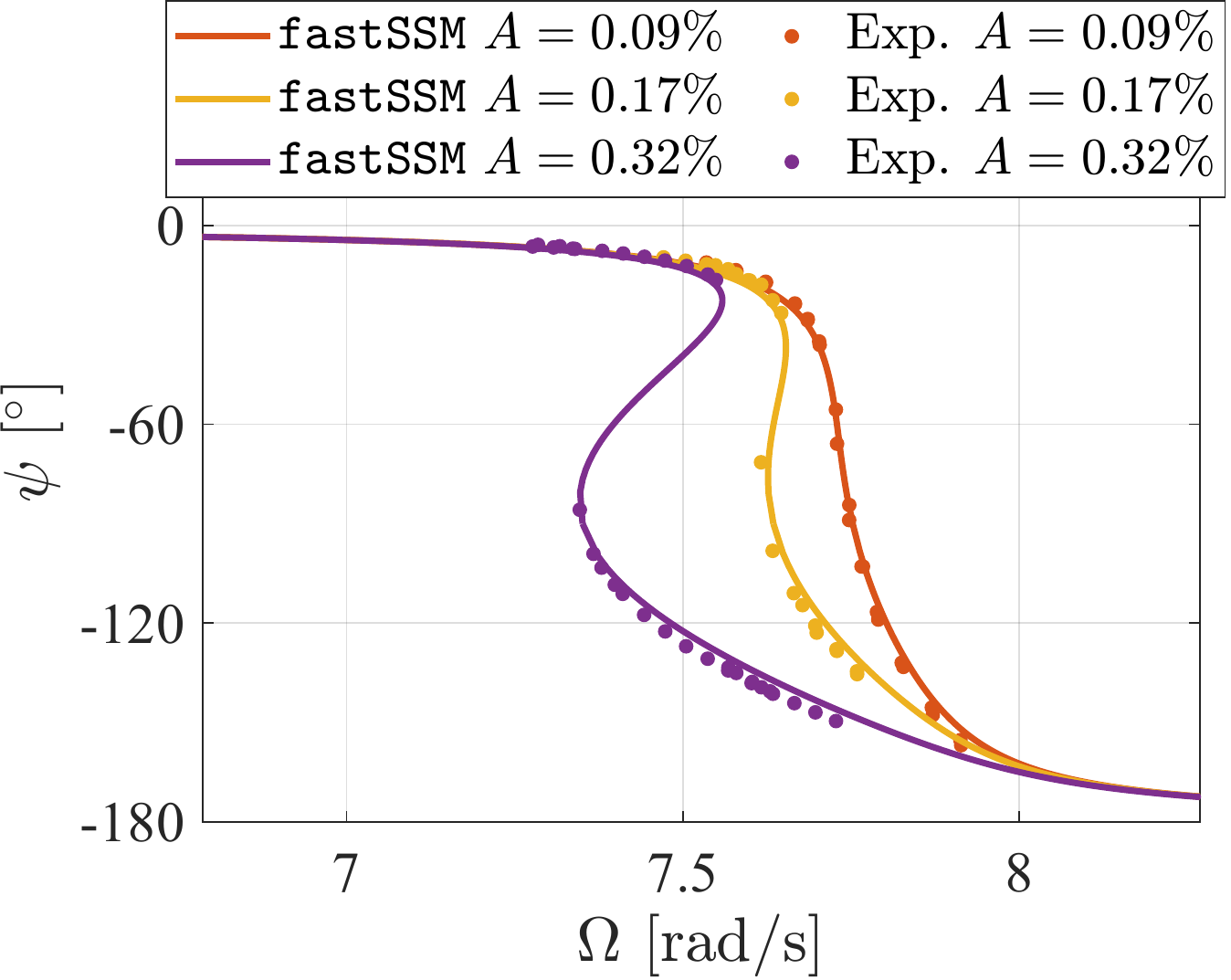}
         \caption{}
         \label{fig:frcphase}
     \end{subfigure}}
     \parbox{0.32\linewidth}{
     \begin{subfigure}{1\linewidth}
         \includegraphics[width=\linewidth]{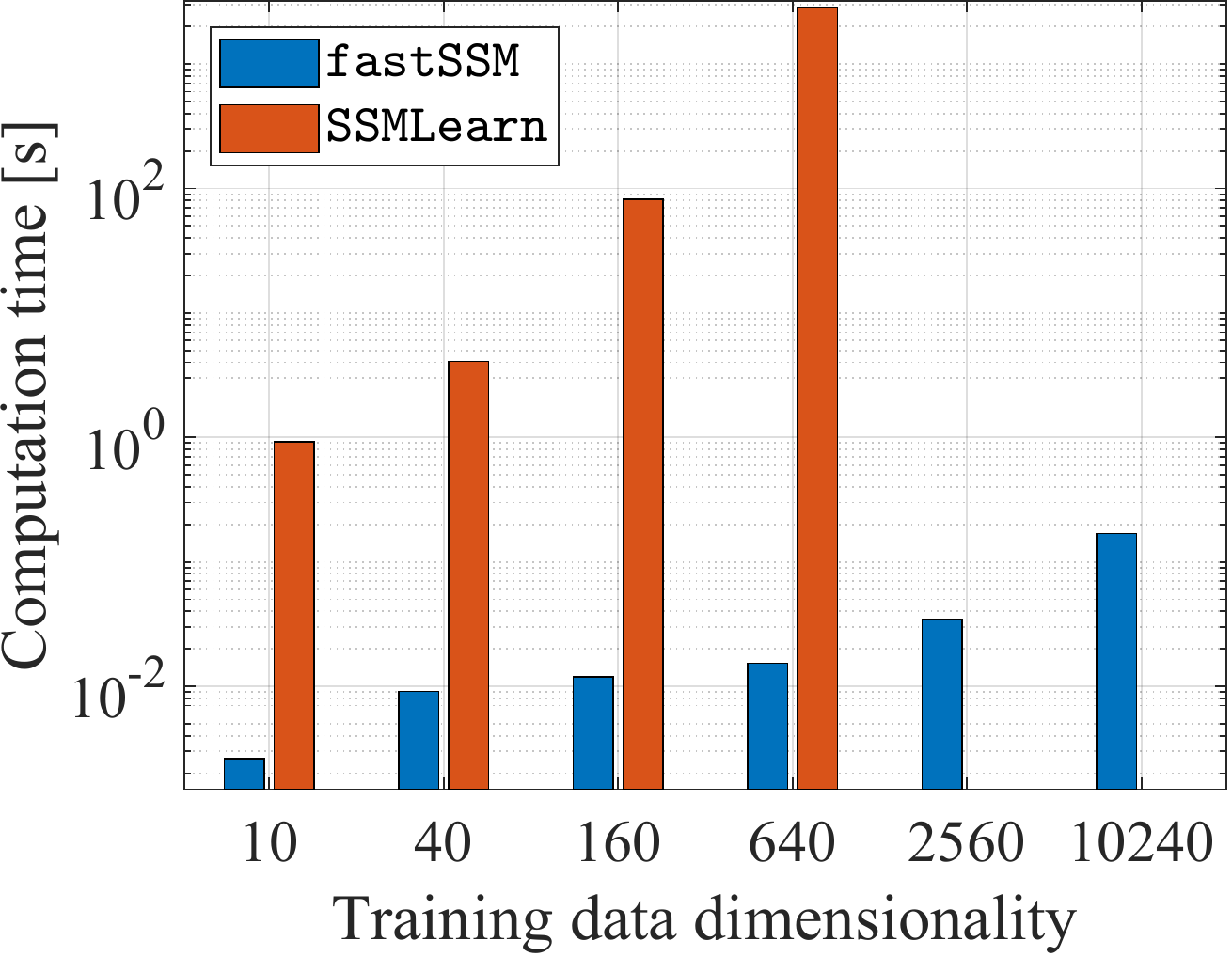}
         \caption{}
         \label{fig:comptime:ndof}
     \end{subfigure}
     \begin{subfigure}{1\linewidth}
         \includegraphics[width=\linewidth]{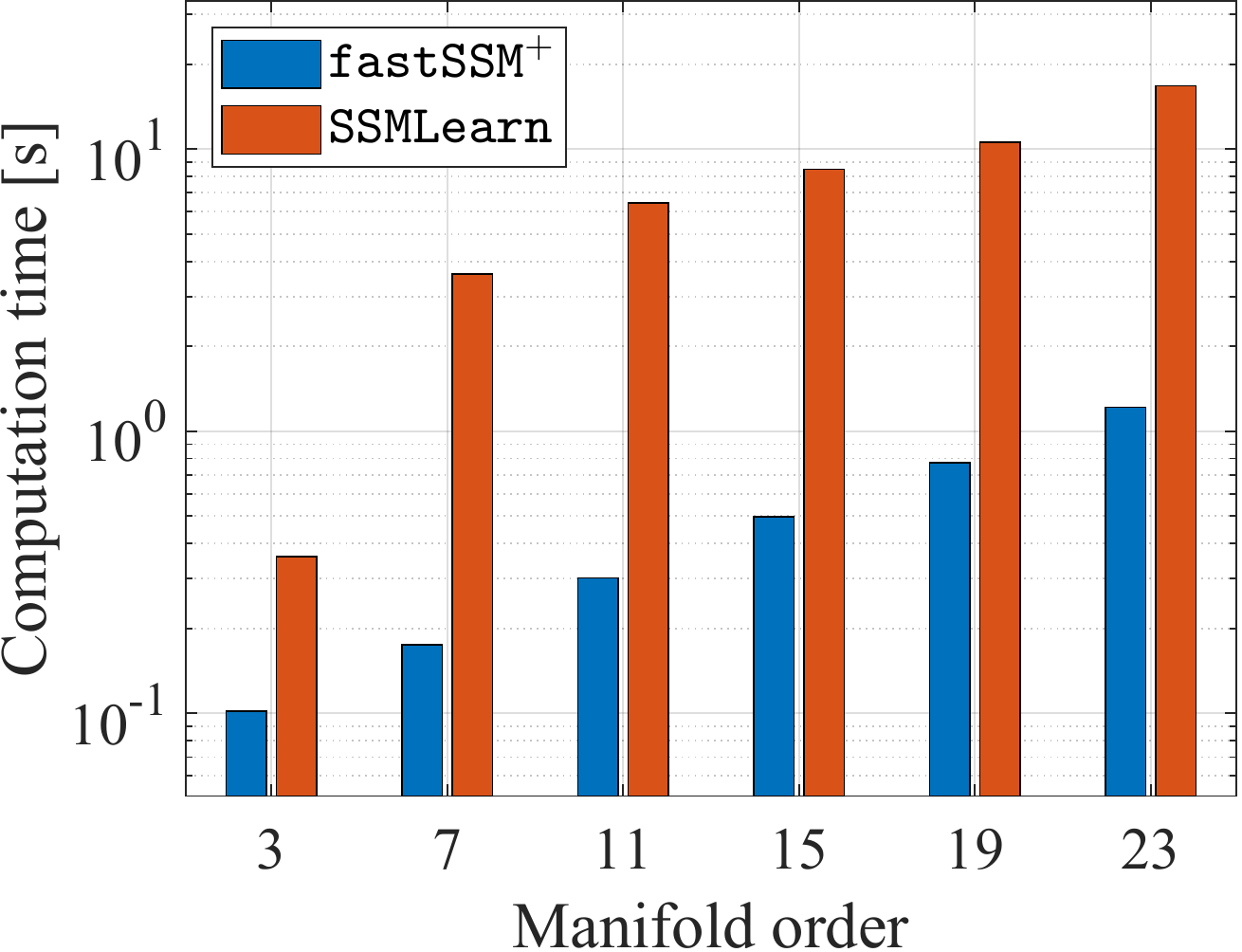}
         \caption{}
         \label{fig:comptime:ordo}
     \end{subfigure}}
     \caption{
     (\subref{fig:ssm}) A visualization of the identified 2D SSM, with the wall surface amplitude on the vertical axis. 
     (\subref{fig:backbones}) The computed instantaneous damping and frequency show significant nonlinearity.
     (\subref{fig:frcamp}) Prediction of the forced response amplitude and (\subref{fig:frcphase}) phase lag compared with experimental results for three different forcing amplitudes.
     (\subref{fig:comptime:ndof}) The computational cost of manifold fitting is orders of magnitude lower with \mSSM{} than with \SSML{}, and scales better with higher dimensionality. 
     (\subref{fig:comptime:ordo}) \mSSMp{} normal form computation is also faster, although \SSML{} tends to require a lower order for convergence.
     }
     \label{fig:sloshing2}
\end{figure*}

In Ref.~\onlinecite{cenedese21}, \SSML{} was applied to experimental sloshing data from Ref.~\onlinecite{bauerlein21}. 
With nonlinear models of both frequency and damping, the forced response of the first oscillation mode was accurately predicted from unforced decaying center of mass signals. 
Here, we will apply \mSSM{} to a much higher-dimensional observable space that includes fine-resolved surface profile data.

The experiments described in Ref.~\onlinecite{bauerlein21} were performed in a rectangular tank of width $w=500$ mm and depth $l=50$ mm, partially filled with water up to a height of 400 mm. 
The tank was mounted on a moving platform excited harmonically by a motor at different frequencies $\Omega$ and amplitudes $A$. 
A camera was mounted on the moving platform, and the surface profile $\vct h$ was detected via image processing with the sampling time $\Delta t = 0.033$ s. 
Figure \ref{fig:experimentsetup} displays the experimental setup. 
%The flow is verified to be quasi-two-dimensional by stereoscopic particle image velocimetry with two high-speed cameras

From the surface profile, the sloshing amplitude can be quantified by computing the horizontal position $\hat{X}$ of the liquid's center of mass at each time, normalized by the tank width. 
$\hat{X}$ is a physically meaningful quantity, relevant for engineering applications and robust against small image evaluation errors. 
The tank was excited at the tested frequencies until a steady state was reached, and the driving was turned off. 
The oscillation amplitude then decayed along the backbone curve defined in (\ref{eq:backbone}).

Here we append the $\hat{X}$ signal with the measurements of the surface elevation $\vct h$ at $1\,531$ points. 
We delay-embed these signals using 10 subsequent measurements to create a $15\,320$-dimensional observable space, in which \mSSM{} identifies a 2-dimensional, 7th-order manifold, shown in \autoref{fig:ssm}. 
We identify the reduced dynamics on the SSM up to 3rd order, and then compute its 3rd-order normal form
\begin{equation}\label{eq:sloshingnf}
\begin{aligned}
	\dot{\rho} &= -0.062 \rho -0.019 \rho^3, \\
	\rho\dot{\theta} &= 7.81\rho -0.628 \rho^3.
\end{aligned}
\end{equation}
The backbone curves obtained from this normal form are shown in \autoref{fig:backbones}. 

In Figures \ref{fig:frcamp} and \ref{fig:frcphase}, we compute forced response curves for the center of mass position, $\hat{X}$, and compare to its experimentally generated values obtained along frequency sweeps from both directions. 
We find the response prediction from \mSSM{} to be accurate.
In particular, the nonlinear damping term of the normal form helps capturing the width of the forced response curve at higher amplitudes. 

\begin{figure*}
    \centering
    \parbox{0.85\linewidth}{
    \parbox{0.373\linewidth}{
    \begin{subfigure}{1\linewidth}
        \includegraphics[width=\linewidth]{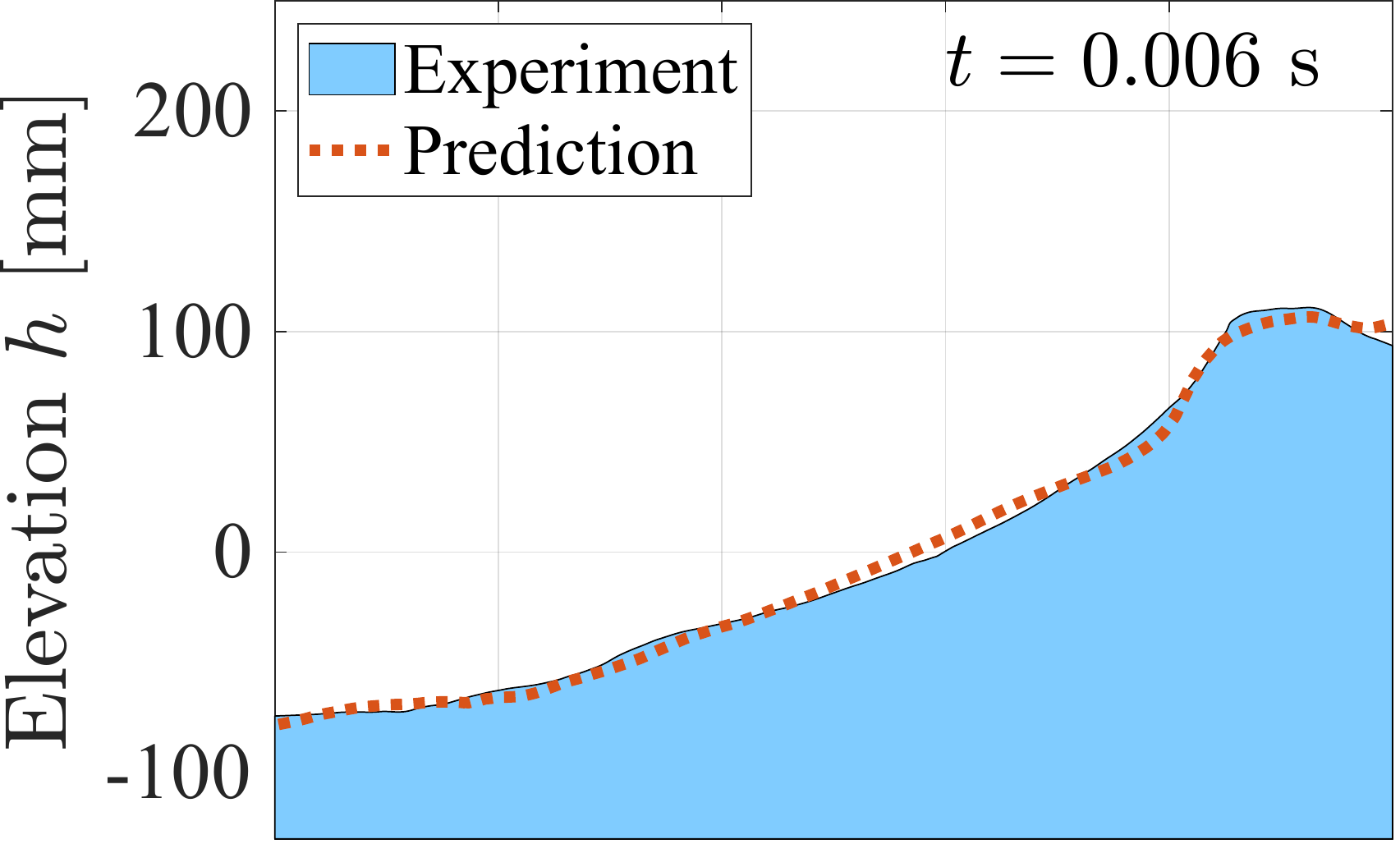}
    \end{subfigure}
    \begin{subfigure}{1\linewidth}
        \includegraphics[width=\linewidth]{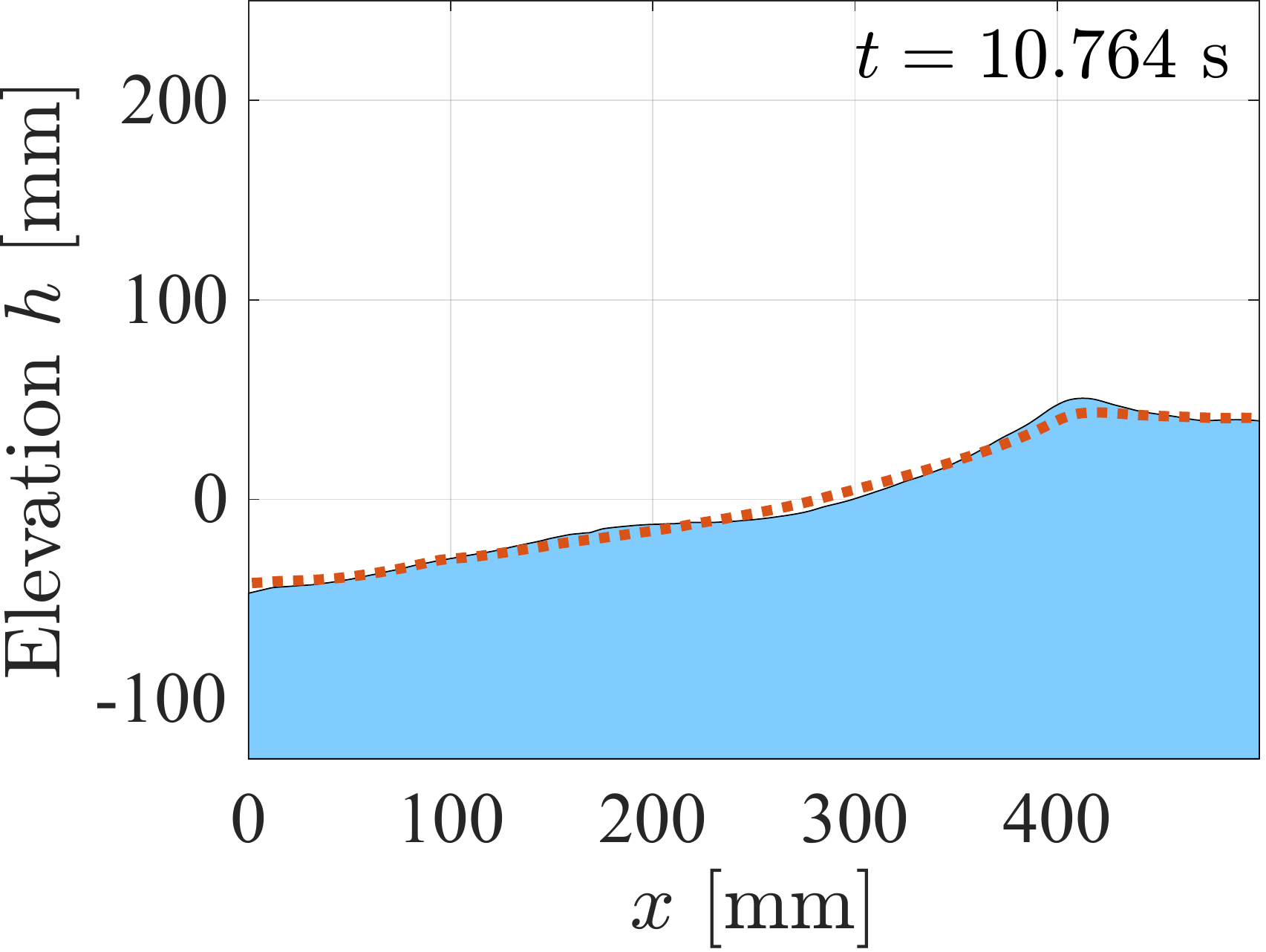}
    \end{subfigure}
	}
    \parbox{0.305\linewidth}{
    \begin{subfigure}{0.98\linewidth}
        \includegraphics[width=\linewidth]{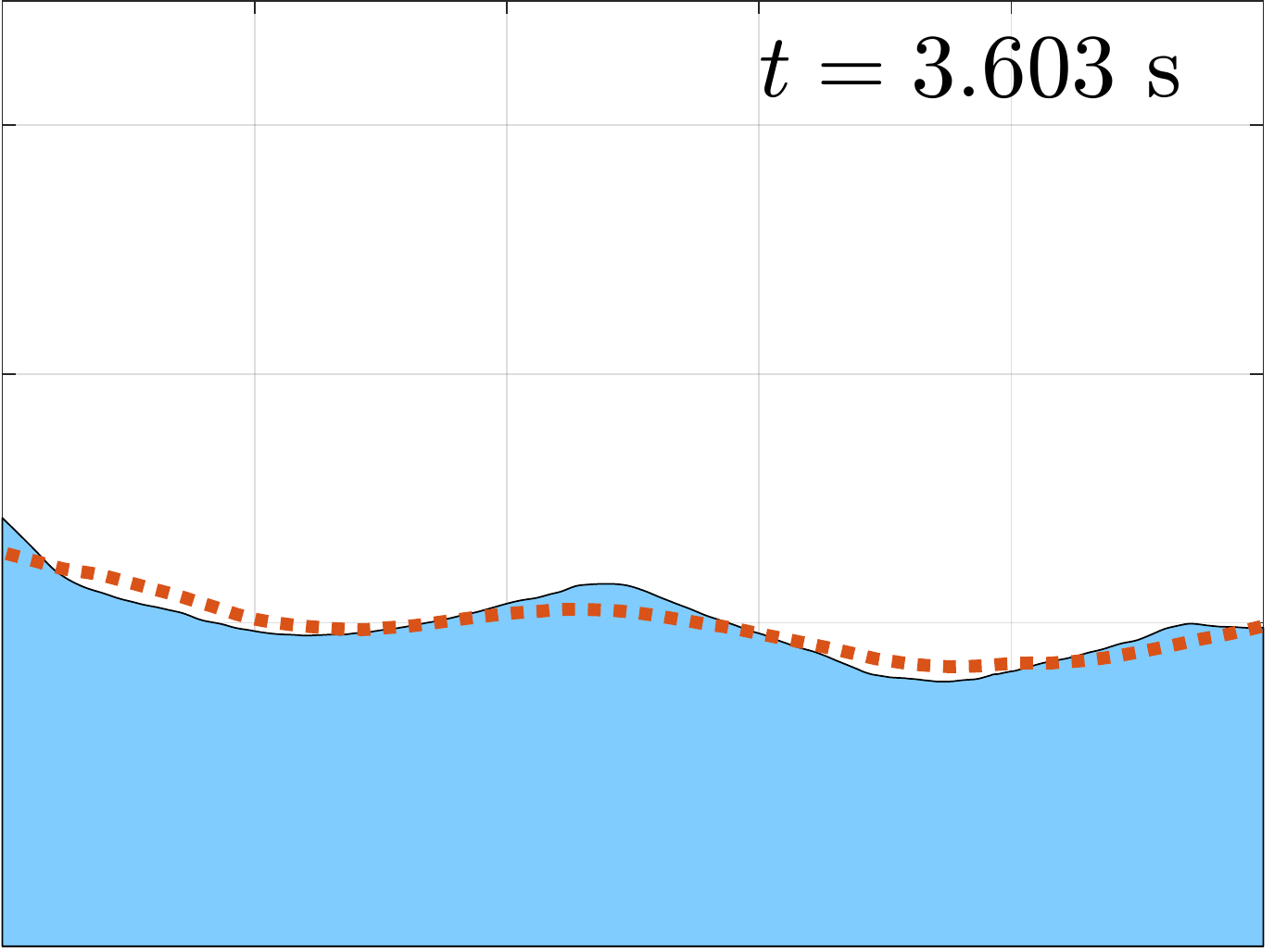}
    \end{subfigure}
    \begin{subfigure}{1\linewidth}
        \includegraphics[width=\linewidth]{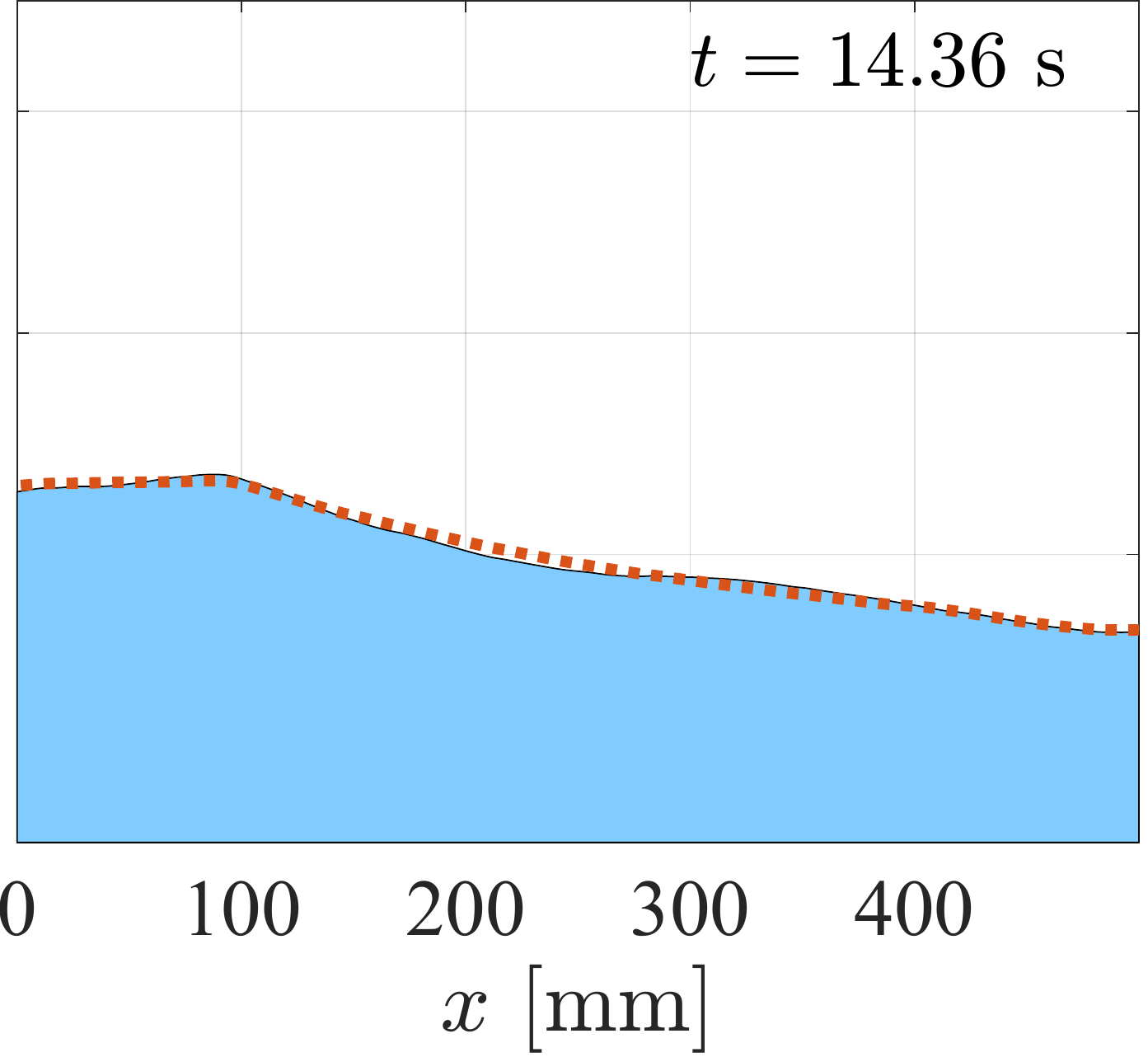}
    \end{subfigure}
    }
    \parbox{0.305\linewidth}{
    \begin{subfigure}{0.98\linewidth}
        \includegraphics[width=\linewidth]{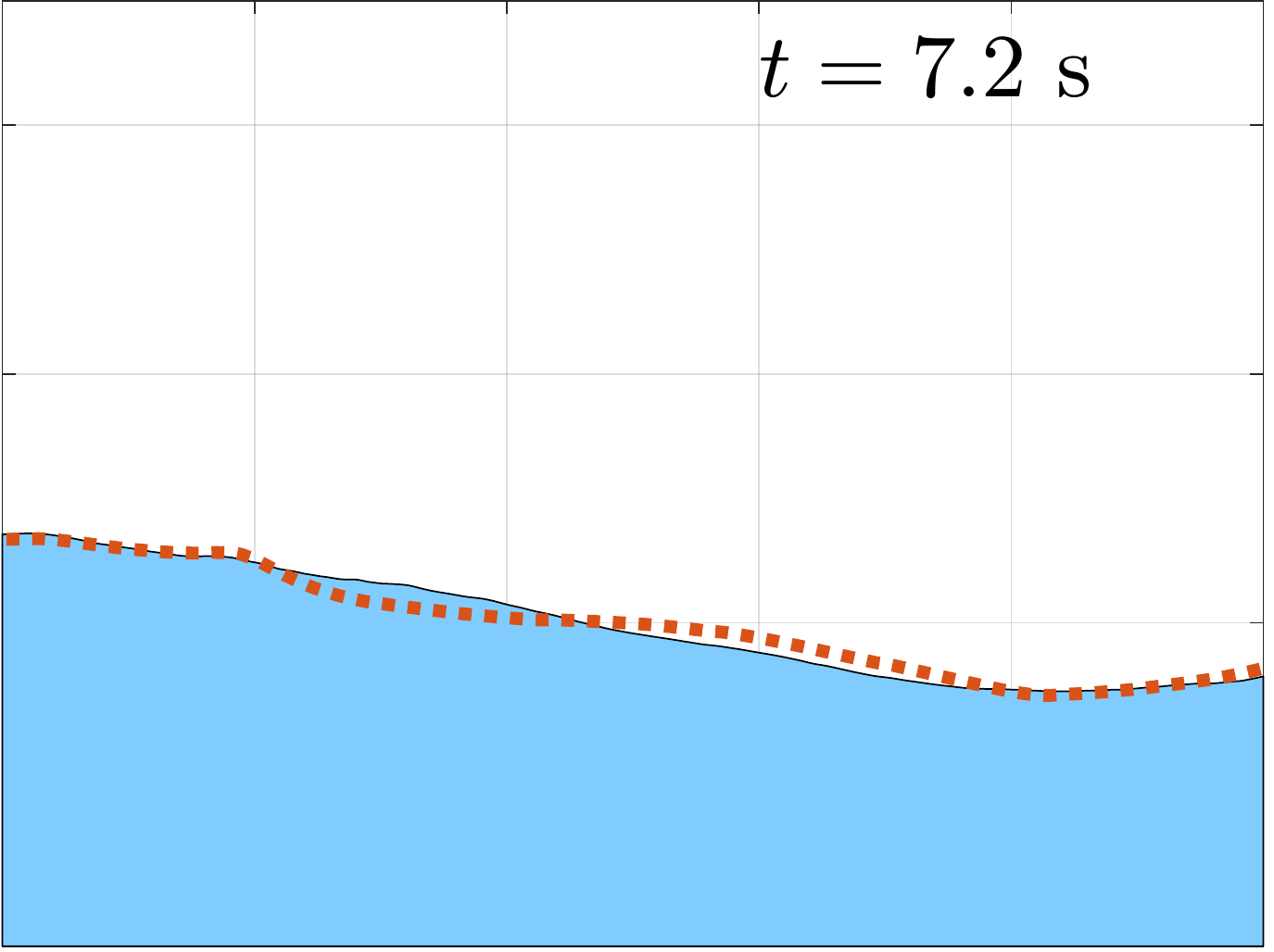}
    \end{subfigure}
    \begin{subfigure}{1\linewidth}
        \includegraphics[width=\linewidth]{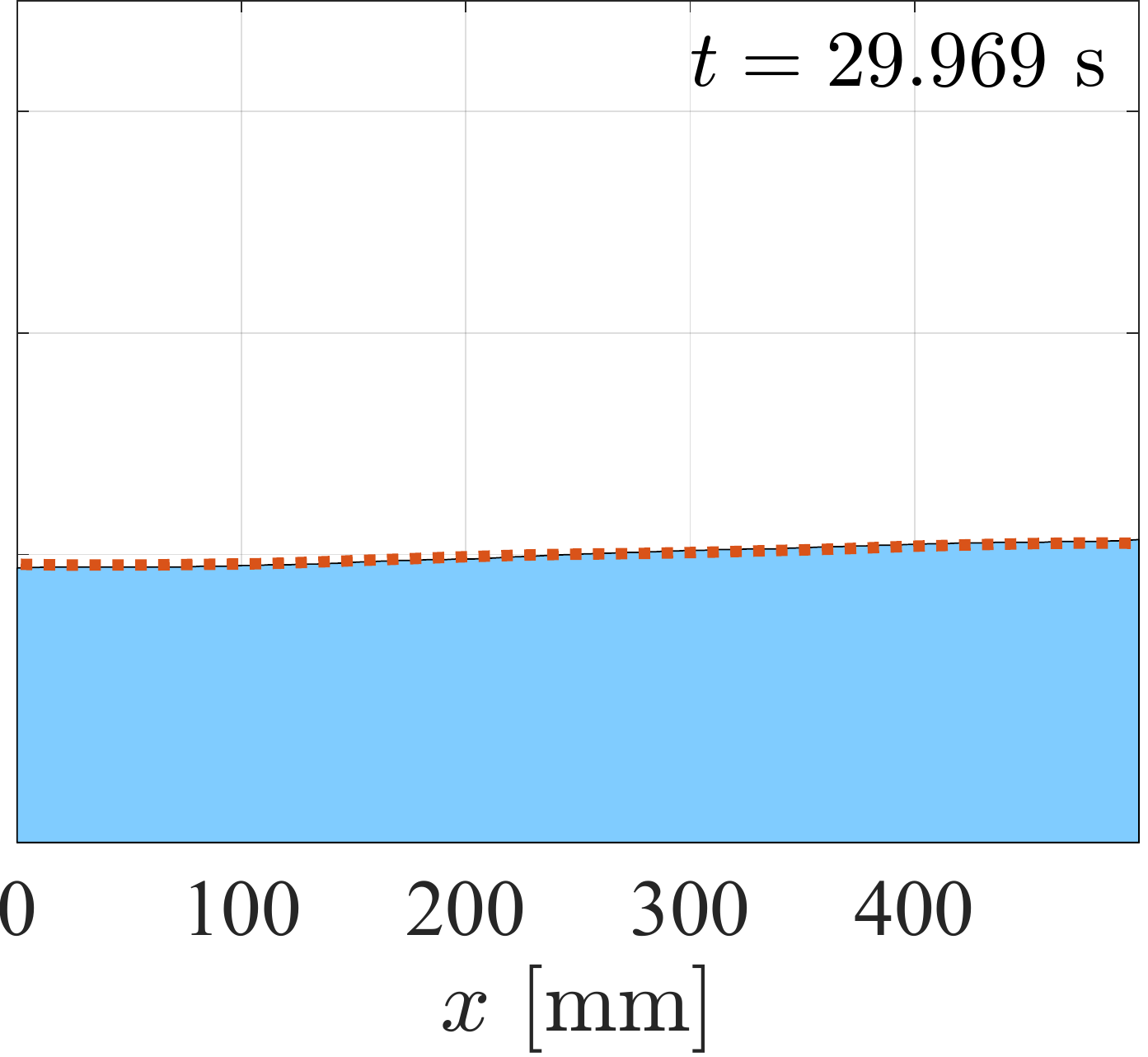}
    \end{subfigure}
    }
    }
    \caption{Snapshots of the \mSSM{} prediction of the surface profile decay. 
    The initial condition is transformed to the SSM-reduced system, which is integrated to produce a 2D trajectory in normal form coordinates. 
    By mapping the trajectory back to the observable space, we can predict the full surface profile and observe close agreement with the experiments in phase, amplitude, and shape.}
    \label{fig:surfaceprofile}
\end{figure*}

We also use the normal form (\ref{eq:sloshingnf}) to predict the development of the decaying full surface profile $\vct h$ in \autoref{fig:surfaceprofile}. 
Here, we take the initial surface profile and transform it to an initial condition in the normal form coordinates.
We integrate the normal form from this initial condition to predict its development in the observable space. 
We observe that the prediction is in close agreement with experiments, yielding a total NMTE of 2.05~\% over the entire phase space.

Finally, we compare execution times for \mSSM{} and \SSML{} when they are trained on the surface profile data. 
These computations were performed on MATLAB version 2020b, installed on an iMac with 2.3 GHz 18-Core Intel Xeon W and 128 GB RAM.

In \autoref{fig:comptime:ndof}, the computational effort of fitting a 2D SSM is plotted against the dimensionality of the training data, i.e., the number of included surface points multiplied by the delay embedding dimension 10. 
Due to its explicit coefficient fitting, \mSSM{} achieves a major speedup, and enables analysis of significantly higher-dimensional data than \SSML{}.

After reduction to the manifold, both \SSML{} and \mSSM{} compute the $\ordo{3}$ normal form in less than a second.
In order to compare the computational effort for higher order normal forms, we apply \mSSMp{}. 
\autoref{fig:comptime:ordo} shows the time required to compute a 2D normal form after fitting the SSM to the sloshing data. 
At a given order, \mSSMp{} is on average 15 times faster than \SSML{}.
While both methods are fast in this example, the difference becomes significant at higher dimensions.
It should be noted, however, that in strongly nonlinear cases, \mSSMp{} may require a higher order normal form to converge, and so the difference to the practitioner may be smaller. 
This difference in the convergence of normal forms is closer examined in the next example.

\subsection{\vk{} beam}\label{subsec:vk}

\begin{figure*}
     \centering
     \parbox{0.5\linewidth}{
     \begin{subfigure}{\linewidth}
         \centering
         \includegraphics[width=0.8\linewidth]{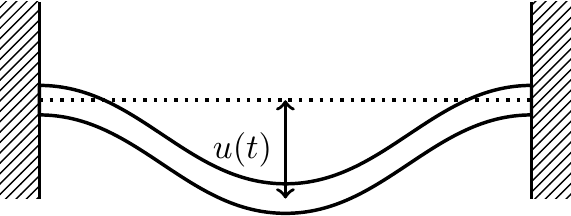}
         \caption{}
         \label{fig:vkmesh}
     \end{subfigure}
     \begin{subfigure}{\linewidth}
         \centering
         \includegraphics[width=\linewidth]{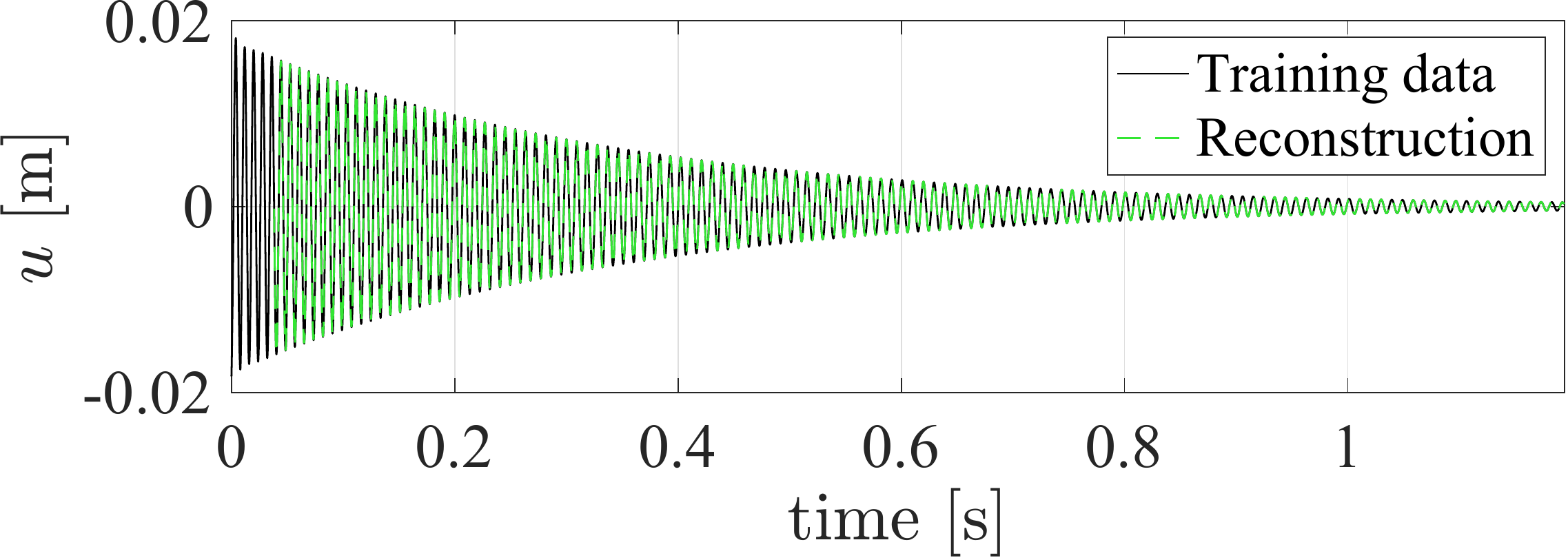}
         \caption{}
         \label{fig:vkdecay}
     \end{subfigure}}
     \begin{subfigure}{0.45\linewidth}
         \includegraphics[width=\linewidth]{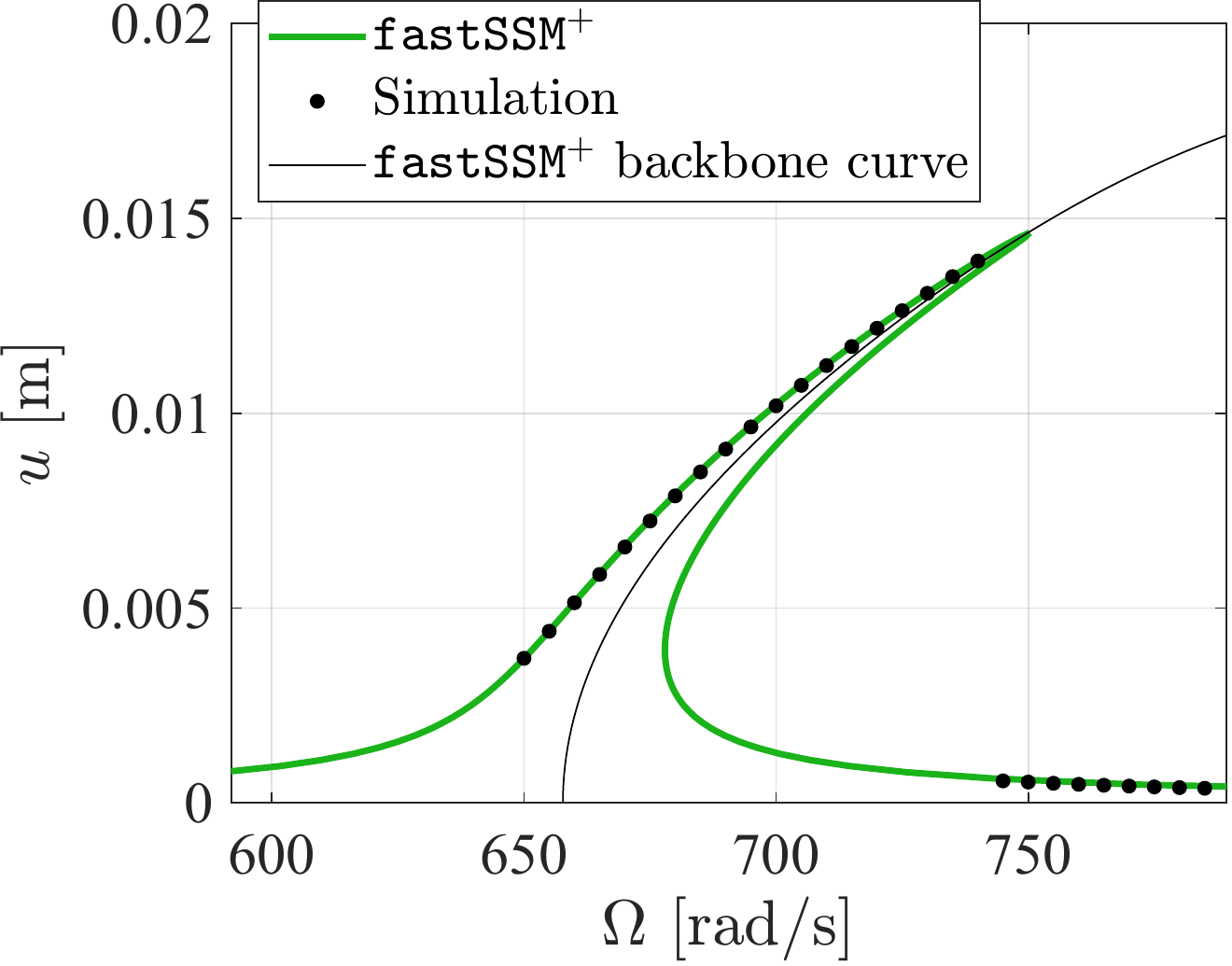}
         \subcaption{}
         \label{fig:mssmvkfrc}
     \end{subfigure}
     \begin{subfigure}{0.32\linewidth}
         \includegraphics[width=\linewidth]{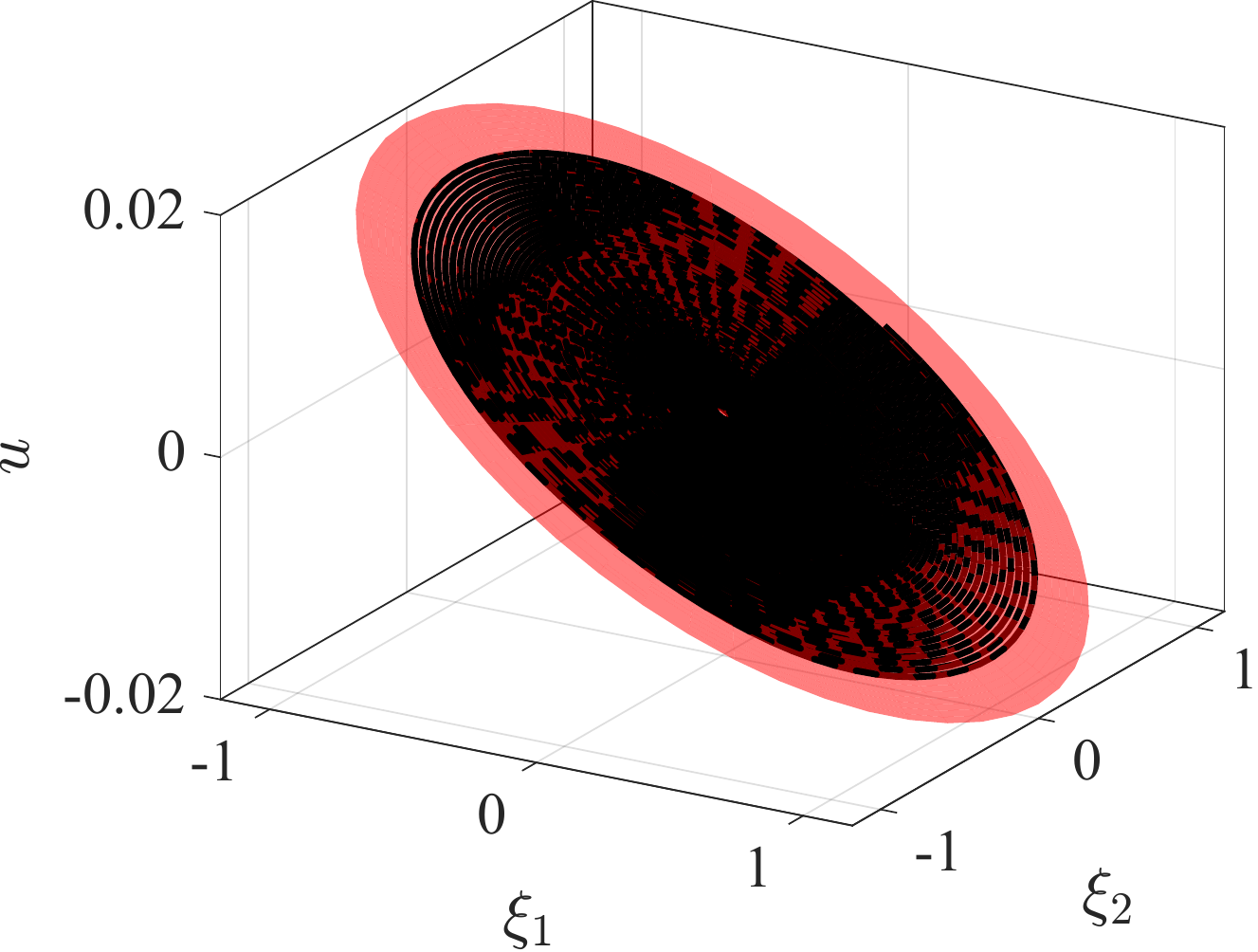}
         \caption{}
         \label{fig:vkssm}
     \end{subfigure}
     \begin{subfigure}{0.32\linewidth}
         \includegraphics[width=\linewidth]{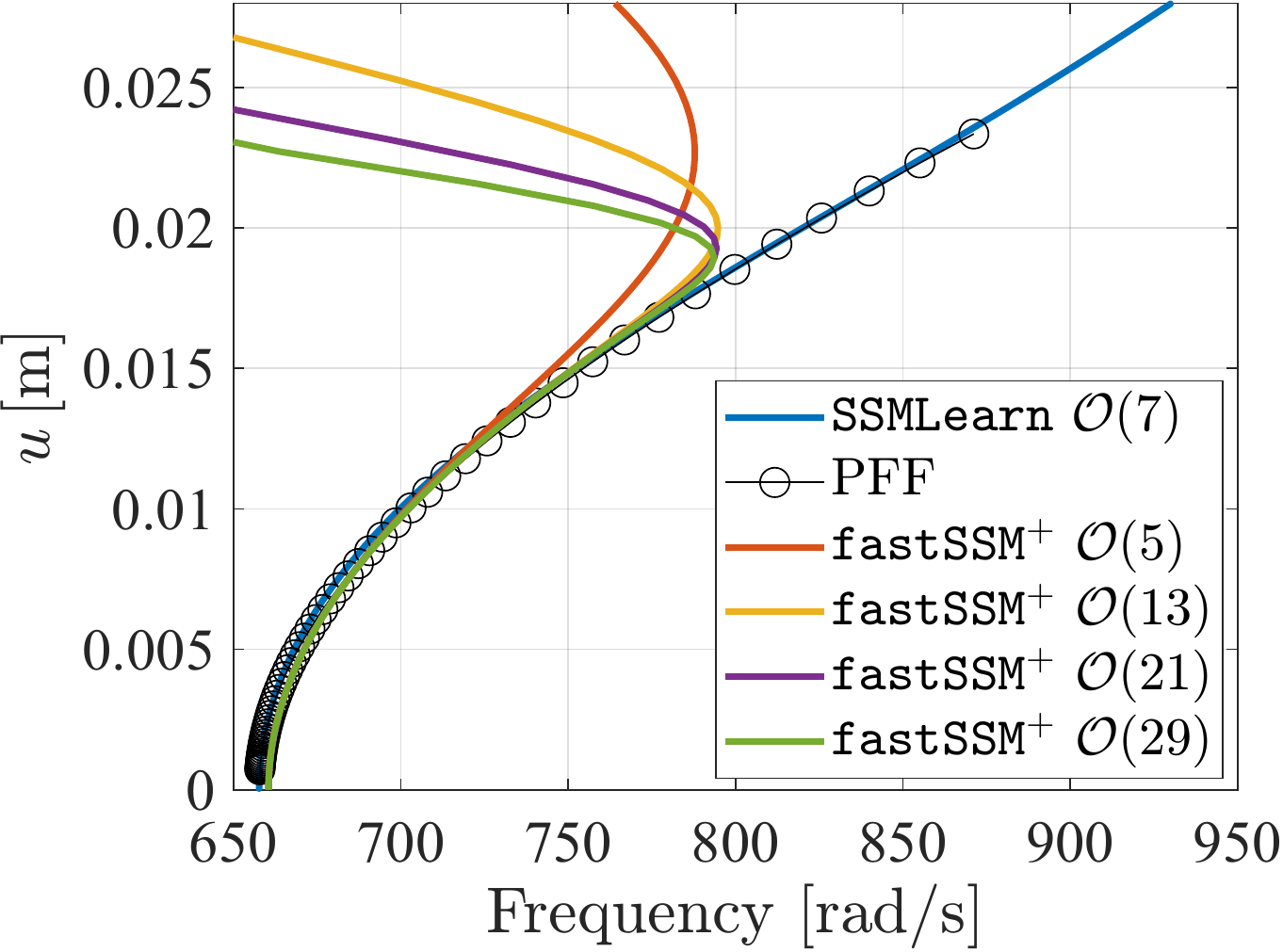}
         \caption{}
         \label{fig:mssmvkconv}
     \end{subfigure}
     \begin{subfigure}{0.32\linewidth}
         \includegraphics[width=\linewidth]{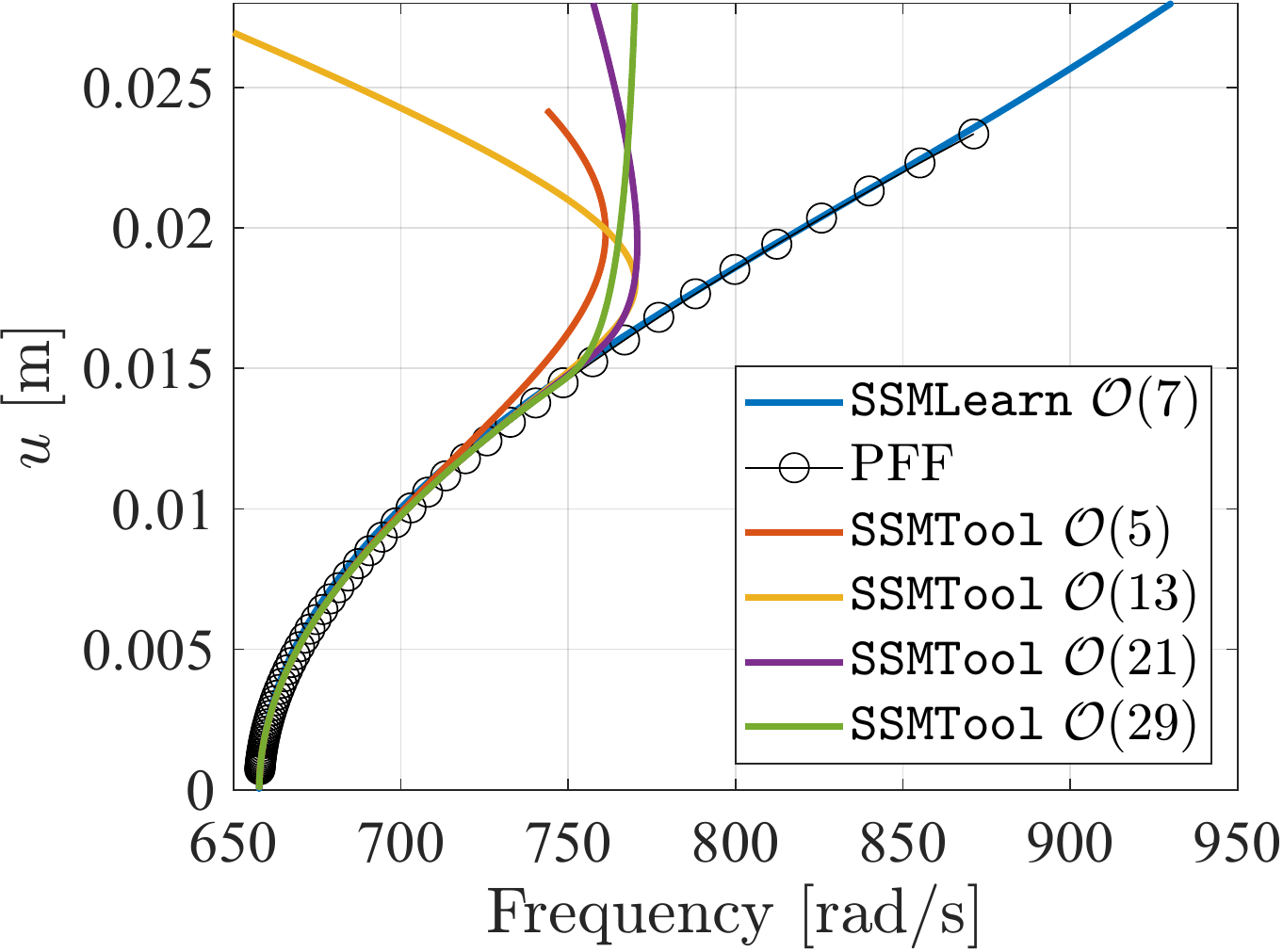}
         \caption{}
         \label{fig:ssmtoolvkconv}
     \end{subfigure}
     \caption{(\subref{fig:vkmesh}) The observable is the midpoint displacement $u(t)$ of a 12-element \vk{} clamped-clamped beam. 
     (\subref{fig:vkdecay}) Based on an $\ordo{5}$ reduced dynamics model, \mSSMp{} returns a normal form that predicts the displacement decay accurately. 
     (\subref{fig:mssmvkfrc}) Prediction of the forced response amplitude.
     (\subref{fig:vkssm}) The 2D SSM visualized in a 3D representation of the observable space.
     Above 20 mm deflection, increased orders of \mSSMp{} (\subref{fig:mssmvkconv}) and \SSMT{} (\subref{fig:ssmtoolvkconv}) no longer improve the backbone curve approximation. Higher amplitudes lie beyond the radius of convergence of the SSM Taylor expansion. \SSML{}, however, continues to approximate the backbone curve well beyond this limit.}
     \label{fig:vonkarman}
\end{figure*}
We consider data from numerical simulations of a finite-element model of a clamped-clamped \vk{} nonlinear beam. \cite{jain18}
Each element has three degrees of freedom: axial deformation $u_0$, transverse deflection $w_0$, and rotation $w_0'$. 
%The displacement field is given by
%\begin{equation}
%    \begin{aligned}
%        u_1(x_1,x_3) &= u_0(x_1) - x_3w_0'(x_1) \\
%        u_2(x_1,x_3) &= 0 \\
%        u_3(x_1,x_3) &= w_0(x_1).
%    \end{aligned}
%\end{equation}
The nonlinear \vk{} axial strain approximation is
\begin{equation}
    \epsilon_{11} = u_0'(x_1) + \frac{1}{2}\left(w_0'(x_1)\right)^2 - zw_0''(x_1),
\end{equation}
where the second term sets this model apart from the linear Euler-Bernoulli beam. 
The axial stress is modeled as
\begin{equation}
    \sigma = E\epsilon_{11} + \kappa \dot{\epsilon}_{11},
\end{equation}
where $E$ is the Young's modulus and $\kappa$ is the material rate of viscous damping.

After a convergence analysis, we set the number of elements to 12, which results in a 33-degree of freedom mechanical system, i.e., a 66-dimensional phase space. 
As initial condition, we compute the response to a static transverse load of 14 kN at the midpoint, removed at the simulation start. 
Our observable function is the midpoint displacement, and the objective is to reconstruct the SSM and its dynamics in the observable space by delay-embedding the signal. 
A sketch of the system is shown in \autoref{fig:vkmesh}.

We set $E=70$ GPa, $\kappa=1.0\times 10^6$ $\mathrm{Pa\cdot s}$, length 1000 mm, width 50 mm, and thickness 20 mm. 
The sampling time is $\Delta t = 0.0955$ s. 
To satisfy the conditions of Takens' embedding theorem, we set the delay embedding dimension to $p=5$. 
The maximum displacement in the training data is 15.9 mm. 

The cubic normal form in \mSSM{} is insufficient to describe the higher amplitude oscillations, so we deploy \mSSMp{} to compute a higher-order normal form.
After training on the generated trajectory, \mSSMp{} outputs a normal form that we use for forced response prediction. 
We use a 1st-order SSM, 5th-order SSM-reduced dynamics, and obtain an 11th-order normal form
\begin{subequations}\label{eq:vknf}
\begin{equation}
\begin{aligned}
	\dot{\rho} &= -3.12 \rho-0.899 \rho^3+4.12 \rho^5-\\&12.4 \rho^7+16.0 \rho^9-17.7 \rho^{11},
\end{aligned}
\end{equation}
\begin{equation}
\begin{aligned}
	\rho\dot{\theta} &= 658 \rho+263 \rho^3-155 \rho^5+\\&151 \rho^7-181 \rho^9+237 \rho^{11},
\end{aligned}
\end{equation}
\end{subequations}
which yields a NMTE of 0.63~\% on the training data (\autoref{fig:vkdecay}). 

We compute the FRC and verify it against a numerical integration frequency sweep in \autoref{fig:mssmvkfrc}. 
Clearly, the autonomous SSM computed by \mSSMp{} can predict the forced response very accurately even with strong nonlinearities. 
\autoref{fig:vkssm} shows a representation of the SSM geometry, which,  as predicted in \autoref{subsec:embed}, is almost flat due to the small timelag.

\begin{figure*}
     \centering
     \begin{subfigure}{0.34\linewidth}
         \includegraphics[width=\linewidth]{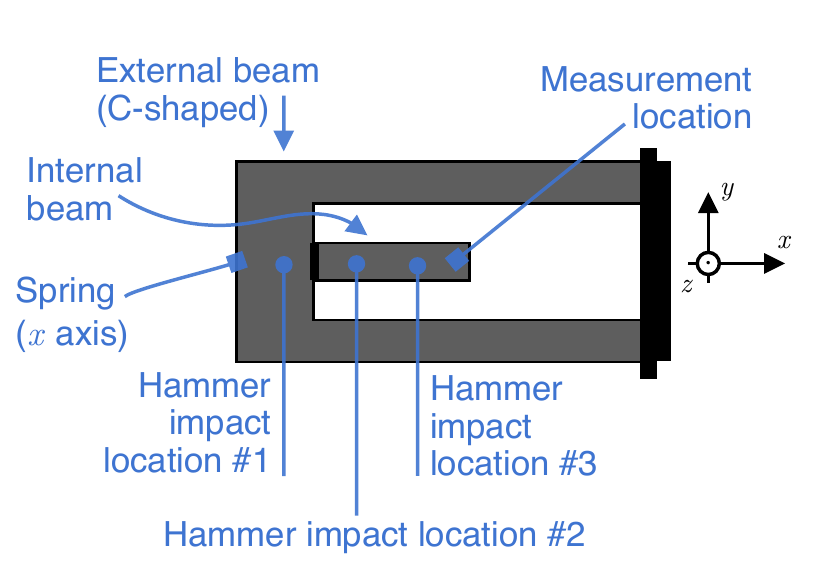}
         \caption{}
         \label{fig:resbeam:sketch}
     \end{subfigure}
     \begin{subfigure}{0.32\linewidth}
         \includegraphics[width=\linewidth]{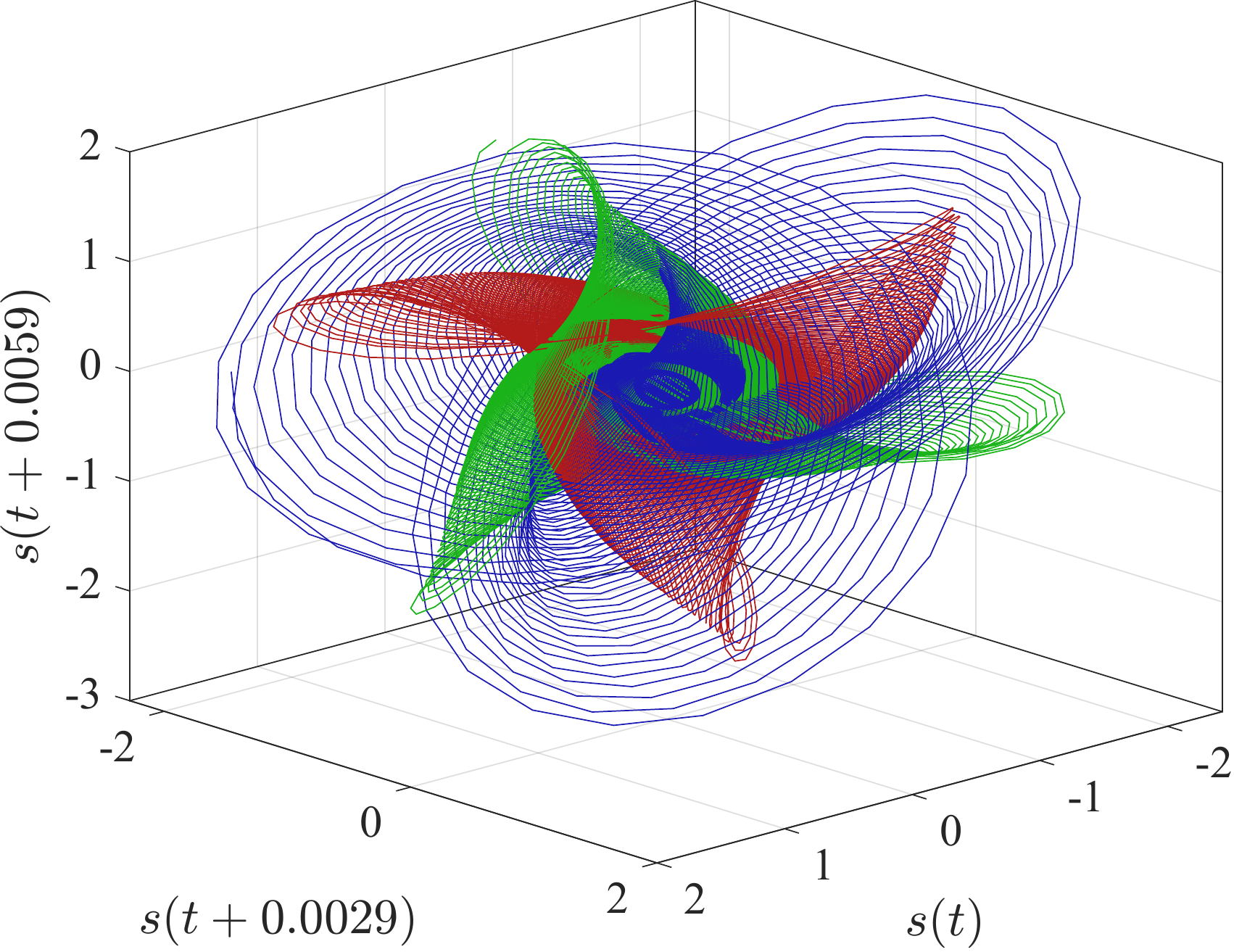}
         \caption{}
         \label{fig:resbeam:y}
     \end{subfigure}
     \begin{subfigure}{0.32\linewidth}
         \includegraphics[width=\linewidth]{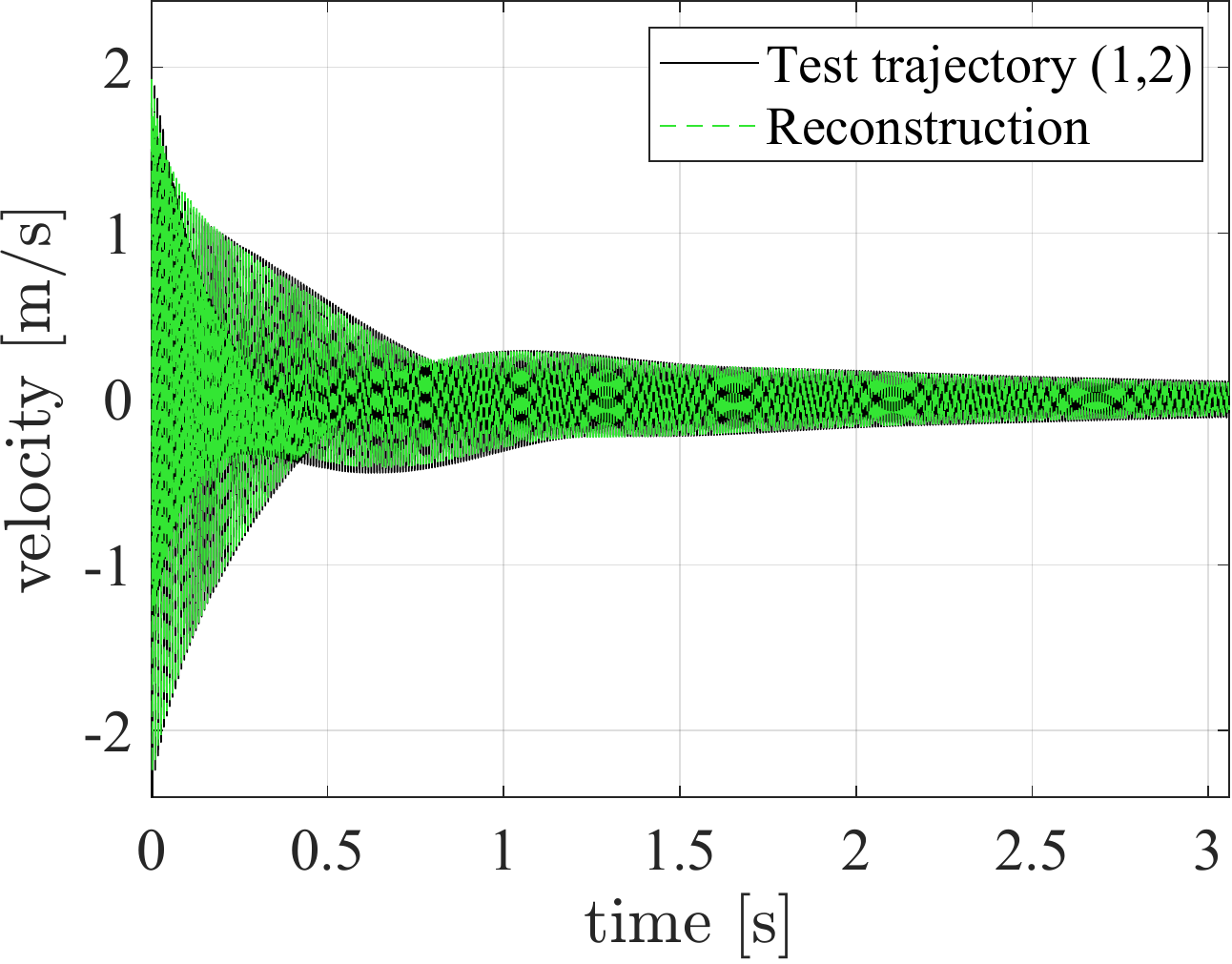}
         \caption{}
         \label{fig:resbeam:rec}
     \end{subfigure}
     \begin{subfigure}{0.34\linewidth}
         \includegraphics[width=\linewidth]{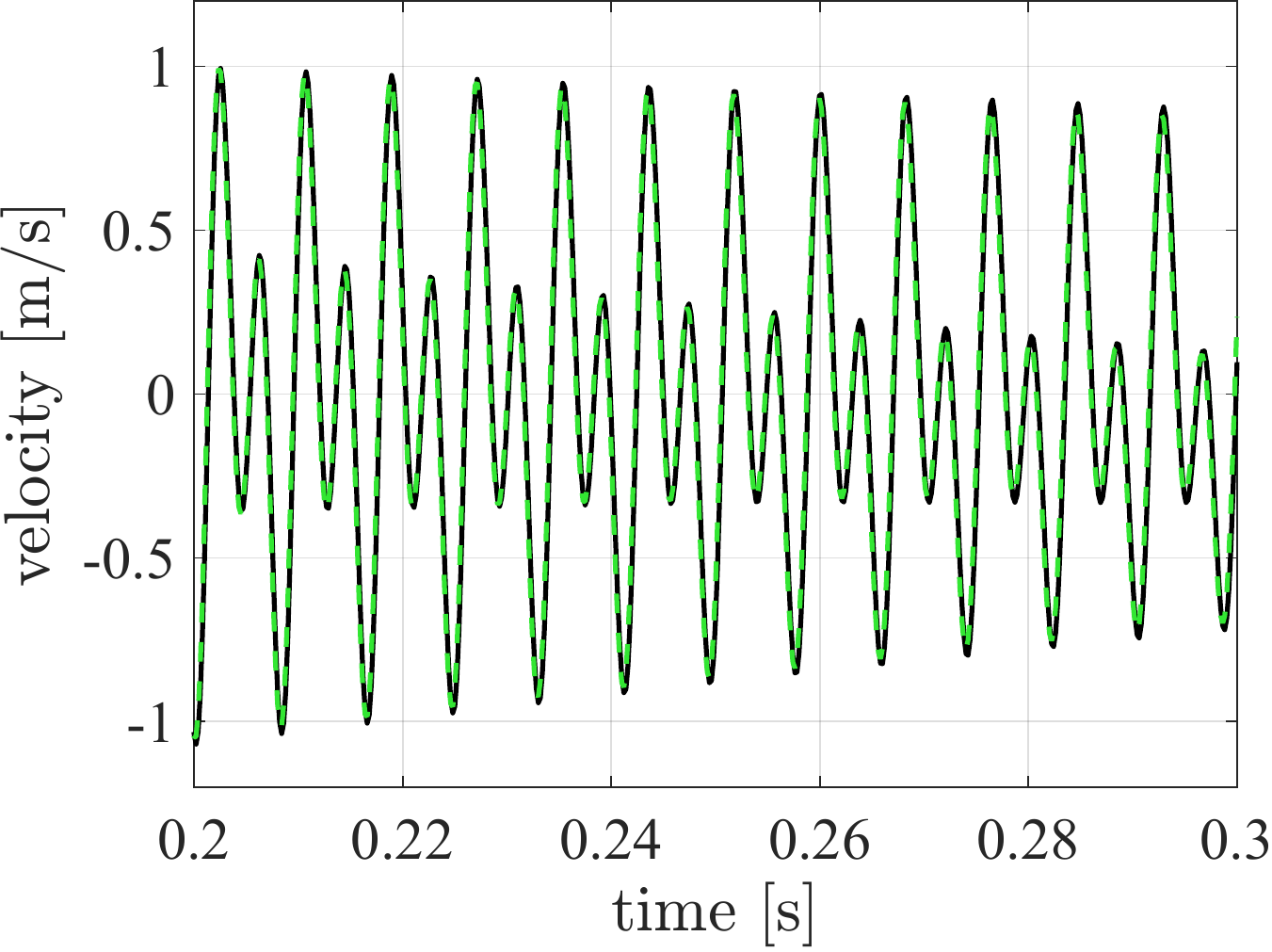}
         \caption{}
         \label{fig:resbeam:reczoom}
     \end{subfigure}
     \begin{subfigure}{0.29\linewidth}
         \includegraphics[width=\linewidth]{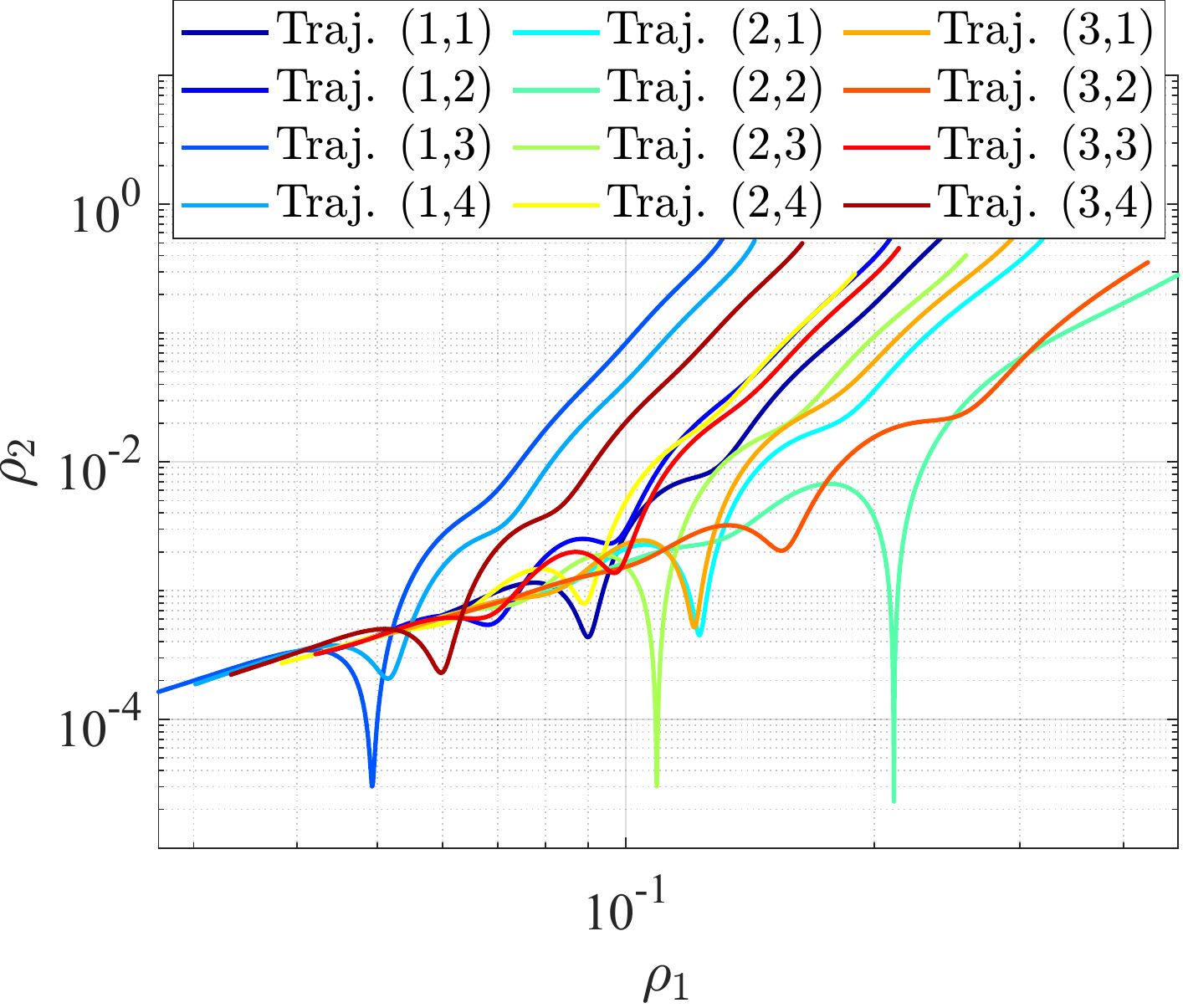}
         \subcaption{}
         \label{fig:resbeam:amps}
     \end{subfigure}
     \begin{subfigure}{0.35\linewidth}
         \includegraphics[width=\linewidth]{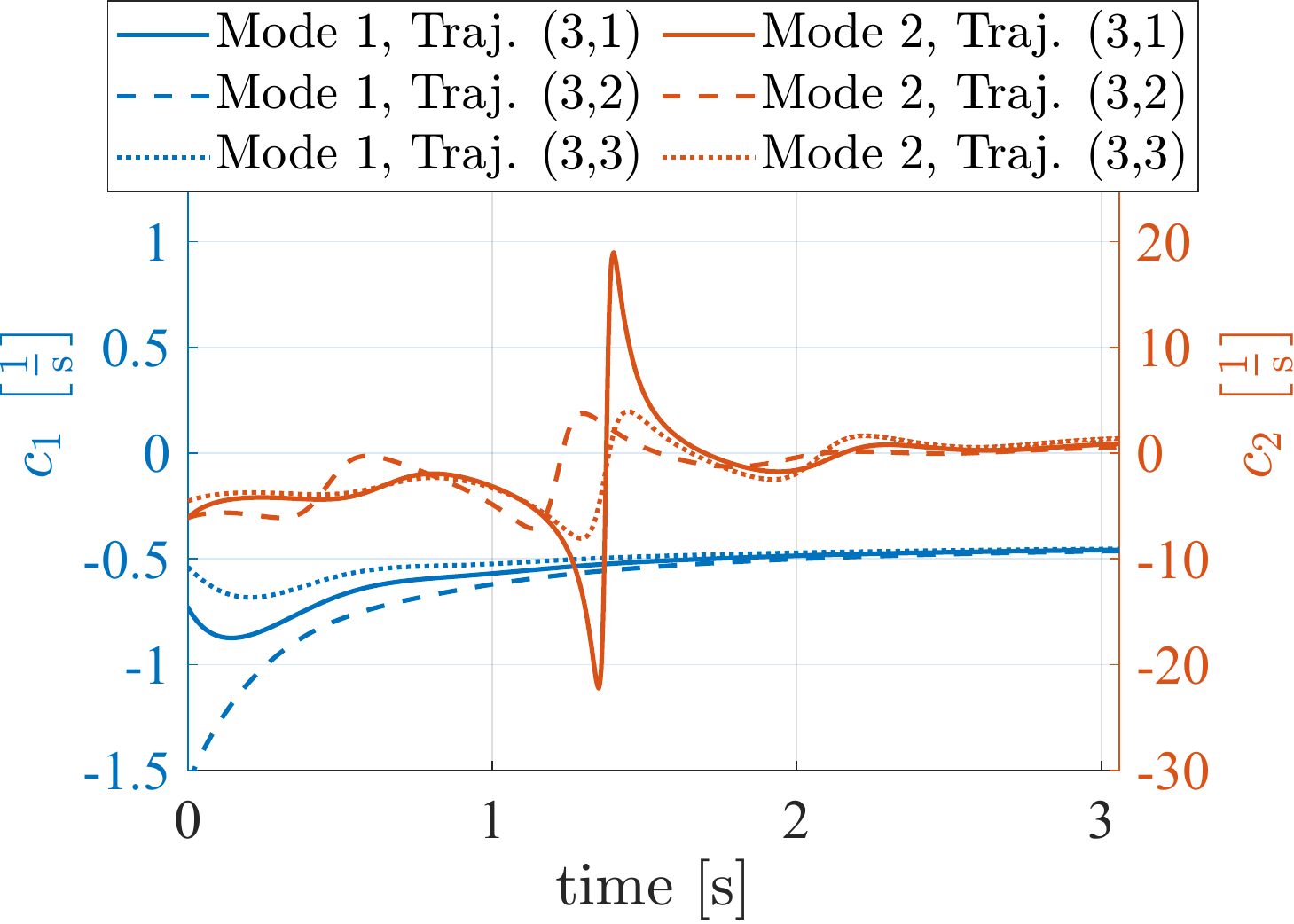}
         \caption{}
         \label{fig:resbeam:damp}
     \end{subfigure}
     \caption{(\subref{fig:resbeam:sketch}) Sketch of the experimental setup for the resonant beam tester. \cite{cenedese21b}
     (\subref{fig:resbeam:y}) The test trajectories embedded in the delay-embedded observable space, in which we fit the SSM.
     (\subref{fig:resbeam:rec}, \subref{fig:resbeam:reczoom}) The reconstruction from the obtained 4D normal form dynamics agrees with the decay of measured velocities.
     (\subref{fig:resbeam:amps}) Plotting the modal amplitudes $\rho$ against one another shows that $\rho_2$ dies out more quickly than $\rho_1$, but, due to modal interaction, the decay rates vary.
     (\subref{fig:resbeam:damp}) The instantaneous damping of the second mode is at times positive, indicating modal energy interchange.}
     \label{fig:resbeam}
\end{figure*}

There is, however, a limit to the range of validity of the normal form (\ref{eq:vknf}), and so the forced response prediction is not guaranteed to be accurate for arbitrarily high amplitudes. 
To explore this limitation, we increase the initial point load to 35 kN, equivalent to a maximum displacement of 24.4 mm in the training data. 
We plot backbone curves computed with \mSSMp{} at increasing orders in \autoref{fig:mssmvkconv}. 
For reference, we also compute an approximation of the instantaneous frequency with Peak Finding and Fitting (PFF) from Ref.~\onlinecite{jin20}.
Above approximately 20 mm transverse displacement of the beam midpoint, increasing the order of the normal form computation no longer improves the model, as the radius of analyticity seems to be reached. 
On the other hand, the \SSML{} model remains valid far beyond this limit. 
This is a clear advantage of the numerical normal form approach in \SSML{} over the analytical computation in \mSSMp{}.

A similar convergence study is shown in \autoref{fig:ssmtoolvkconv} for increasing orders of \SSMT{} computation on the full finite-element system of equations. 
The \SSML{} backbone curve agrees with the PFF estimate at a larger range than \SSMT{} converges. 
This shows that data-driven reduced-order modeling can surpass analytical methods in terms of range of validity, even when the full system is known.

\subsection{Resonant beam experiments}\label{subsec:resbeam} 

Our final dataset comprises velocity measurements from an internally resonant beam structure detailed in Ref.~\onlinecite{cenedese21b}. 
It consists of an internal beam jointed with three bolts to the left midpoint of an external C-shaped beam, which is clamped to ground at its rightmost edges (\autoref{fig:resbeam:sketch}). 
The bolts give rise to nonlinear frictional slip, \cite{eriten13, brake18} and the beam lengths are tuned so that the system has a 1:2 resonance between its slowest out-of-plane bending eigenfrequencies, measured to be 122.4 Hz and 243.4 Hz. 
Vibrations are induced with an impulse hammer at three different impact locations and the transient out-of-plane velocity of the inner beam tip is measured at 5120 Hz. 
Varying the impact location and force for data diversity, 12 different trajectories were obtained. 
Frequency analysis shows that only the two slowest eigenfrequencies are present in the signals, so we aim to identify a 4D internally resonant SSM in an observable space for our model reduction.

With more than one mode present in the data, we must judiciously choose the delay embedding parameters.
In the experiment, the sampling time is $\Delta t = 0.0001953$ s, which with $\omega_2=2\pi\times243.4$ rad/s produces $\tau = 5\Delta t$ according to the method in \autoref{subsec:embed}.
This timelag results in a good reconstruction already at the minimum Takens dimension, $p=9$. 
However, increasing the dimensionality further increases modal orthogonality and thus improves trajectory reconstruction. 
Motivated by this, we select a delay embedding dimension of $p=12$, as the reconstructions improve only marginally beyond this number.
We show a representation of the test trajectories in this observable space in \autoref{fig:resbeam:y}.

We select three trajectories as the test set -- one for each hammer impact location -- and use the remaining nine trajectories for training. 
We set the order of expansion for the SSM $m=3$, the SSM-reduced dynamics order $r=3$, and obtain the 3rd-order normal form. 
The results do depend on the order of reduced dynamics, but are insensitive to the orders of computation for the SSM and the normal form on it. 
\mSSMp{} automatically detects an inner resonance from the data and returns the normal form
%\begin{widetext}
\begin{subequations}
\begin{equation}
\begin{aligned}
    \dot{\rho}_{1} &= -0.4426\rho_{1}+0.0837\rho_{1}\rho_{2}^2-5.550\rho_{1}^3+\\
    &\left(0.2683\cos\psi-0.2639\sin\psi\right)\rho_{1}\rho_{2},
\end{aligned}
\end{equation}
\begin{equation}
\begin{aligned}
    \rho_{1}\dot{\theta}_{1} &= 768.9\rho_{1}-0.7541\rho_{1}\rho_{2}^2-17.90\rho_{1}^3+\\
    &\left(-0.2639\cos\psi-0.2683\sin\psi\right)\rho_{1}\rho_{2},
\end{aligned}
\end{equation}
\begin{equation}
\begin{aligned}
    \dot{\rho}_{2} &= -3.125\rho_{2}-4.750\rho_{2}^3-12.00\rho_{1}^2\rho_{2}+\\
    &\left(1.340\cos\psi-0.1161\sin\psi\right)\rho_{1}^2,
\end{aligned}
\end{equation}
\begin{equation}
\begin{aligned}
    \rho_{2}\dot{\theta}_{2} &= 1529\rho_{2}-10.54\rho_{2}^3-13.86\rho_{1}^2\rho_{2}+
    \\&\left(0.1161\cos\psi+1.340\sin\psi\right)\rho_{1}^2,
\end{aligned}
\end{equation}
\end{subequations}
%\end{widetext}
where $\psi=2\theta_{1}-\theta_{2}$. 
The reconstruction of the first test trajectory obtained by integrating the normal form is shown in \autoref{fig:resbeam:rec}, with a zoomed-in version in \autoref{fig:resbeam:reczoom}. 
The NMTE on the test set computed in the first coordinate is 1.28~\%. 

\autoref{fig:resbeam:amps} plots the first modal amplitude, $\rho_1$, against the second one, $\rho_2$, for each simulated trajectory. 
The second mode clearly decays faster and trajectories are eventually dominated by the slower mode $\rho_1$. 
However, due to the modal coupling terms in the normal form, the decay is not monotonic. 
Rather, there is a wealth of different decaying patterns depending on initial conditions.

Finally, \autoref{fig:resbeam:damp} shows the instantaneous damping $\dot{\rho}_1/\rho_1$ and $\dot{\rho}_2/\rho_2$ for the reconstructions corresponding to the third impact location. 
For the slow mode, we observe strong nonlinearity as the instantaneous damping varies dramatically. 
The fluctuations are even larger for the fast mode. 
At times, its instantaneous damping reaches positive values, which indicates that energy is transferred from the slow to the fast mode. \cite{sapsis12}
Overall, the model obtained with \mSSMp{} agrees with the one returned by \SSML{} in Ref.~\onlinecite{cenedese21b}.

\section{Conclusions}\label{sec:conclusions}
We have introduced a fast alternative to a recent data-driven reduction method for nonlinear dynamical systems.
We have also discussed a simplified implementation for the model reduction setting arising most frequently in practice: a two-dimensional slow SSM with underdamped oscillations modelled up to cubic nonlinearities.

Our approach is fundamentally based on the \SSML{} algorithm, \cite{cenedese21, cenedese21b} but turns the fitting of an SSM into an explicitly formulated problem by assuming that its tangent space can be obtained by singular value decomposition on the data. 
Furthermore, we compute the normal form on the manifold explicitly and recursively, rather than solving an implicit minimization problem. 

We have applied this simplified model-reduction algorithm to three datasets: tank sloshing experiments, a nonlinear beam finite-element simulation, and experiments from an internally resonant mechanical structure. 
In all three problems, we obtained a model that accurately predicted the decay of the autonomous system. 
In addition, we have demonstrated that a forcing term can be added to the autonomous model for highly accurate prediction of the forced response amplitude and phase. 
Training on the beam experimental data, our method automatically detected the internal resonance and returned a model that could predict energy repartition among the modes.

The assumptions made for our new algorithms drastically speed up model identification on SSMs from data, increase the possible dimensionality of observable spaces we can tackle, and significantly simplify the code. 
In comparison to the previous method, however, we sacrifice some model accuracy. 
Perhaps more significantly, we have found large differences in normal form convergence domains to the benefit of \SSML{}.

Specifically, we demonstrate on our simulated beam problem that the numerical normal form has a considerably larger range of validity than the analytical normal forms of the full system and consequently those of our new algorithm. 
In other words, data-driven analysis can outperform analytical methods in terms of model validity, even when the full equations of the system are known. 
This suggests that a plausible approach to obtaining a reduced-order model of a finite-element structure would be to simulate the system and train on a small number of trajectories, rather than formulating the full equation system and computing SSMs analytically. 
This idea is a subject of our ongoing work.

\begin{acknowledgments}
We are grateful to Kerstin Avila and Bastian Bäuerlein (U. Bremen) for making their experimental surface profile data from Ref.~\onlinecite{bauerlein21} available to us. We also wish to thank Melih Eriten (U. Wisconsin) for supplying the resonant beam experimental data from Ref.~\onlinecite{cenedese21b}.
\end{acknowledgments}

\section*{Author declarations}
\subsection*{Conflict of interest}
The authors have no conflicts to disclose.

\ifpreprint

\else
\section*{Data availability}
The data that support the findings of this study are openly available in , reference number [reference number].

\clearpage
\begin{widetext}
\appendix*
\section{Code}
\subsection{\mSSM{}}\label{appendix:minissm}
\mSSM{} requires no external packages and can be applied to data out-of-the-box.
For strong nonlinearities or higher-dimensional SSMs, we recommend installing \SSMT{} and running \mSSMp{}.
\mSSM{} is also available at .
\lstinputlisting[label=code:mssm, numbers=none]{fastSSM.m}
%\newpage
\subsection{\mSSMp{}}\label{appendix:minissmplus}
\mSSMp{} is also available at .
Requires installation of \SSMT{}. \cite{SSMTool}
\lstinputlisting[label=code:mssmp, numbers=none]{fastSSMplus.m}
\end{widetext}

\newpage
\fi

\bibliography{SSM_bibliography}% Produces the bibliography via BibTeX.

\end{document}

%% file: packages.tex
\usepackage{hyperref}
\usepackage{graphicx}
\usepackage{subcaption}
\usepackage{amsmath}
\usepackage{amsfonts}
\usepackage{bm}

\usepackage[numbered,framed]{matlab-prettifier}
\usepackage[T1]{fontenc}
\usepackage[ddmmyyyy]{datetime}

%% file: commands.tex
\newcommand{\SSML}{SSMLearn}
\newcommand{\mSSM}{fastSSM}
\newcommand{\mSSMp}{fastSSM$^+$}
\newcommand{\SSMT}{SSMTool}

\newcommand{\real}{\mathrm{Re\, }}
\newcommand{\imag}{\mathrm{Im\, }}
\newcommand{\lbar}{\bar{\lambda}}
\newcommand{\zbar}{\bar{z}}
\newcommand{\N}{\mathbb{N}}
\newcommand{\R}{\mathbb{R}}
\newcommand{\Ce}{\mathbb{C}}
\newcommand{\e}{\mathrm{e}}
\newcommand{\damp}{c}
\newcommand{\freq}{\omega}
\newcommand{\mfd}{\mathcal{M}}
\newcommand{\specsub}{\mathcal{V}_d}
\newcommand{\ordo}[1]{\mathcal{O}(#1)}
\newcommand{\vk}{von Kármán}
\newcommand{\vct}[1]{\bm{#1}}
\renewcommand{\cal}{\mathrm{cal}}